\documentclass[10pt,reqno]{amsart}

\usepackage{amsmath}
\usepackage{amssymb}
\usepackage{amsthm}
\usepackage{graphicx}
\usepackage{psfrag}
\usepackage{mathrsfs}
\usepackage{dsfont}
\usepackage{enumerate}
\usepackage{lscape}

\newcommand{\Nb}{{\mathbb N}}
\newcommand{\Rb}{{\mathbb R}}

\renewcommand{\leq}{\leqslant}
\renewcommand{\geq}{\geqslant}

\providecommand{\bigeq}[1]{\overset{#1}{=\joinrel=}}

\providecommand{\bbigeq}[1]{\overset{#1}{=\joinrel=\joinrel=\joinrel=}}

\DeclareMathOperator{\sgn}{sgn}

\makeatletter
\renewcommand*\env@matrix[1][\arraystretch]{%
  \edef\arraystretch{#1}%
  \hskip -\arraycolsep
  \let\@ifnextchar\new@ifnextchar
  \array{*\c@MaxMatrixCols c}}
\makeatother

\providecommand{\abs}[1]{\lvert#1\rvert} 

\newenvironment{Proof}{\removelastskip \vskip12pt plus 1pt \noindent
{\em Proof.\/}\rm }{\hfill$\square$ \vskip12pt plus 1pt}
\newenvironment{Proofof}{\removelastskip \vskip12pt plus 1pt \noindent
{\em Proof\/}\rm }{\hfill$\square$ \vskip12pt plus 1pt}

\newtheorem{theorem}{Theorem}[section]
\newtheorem{lemma}[theorem]{Lemma}
\newtheorem{prop}[theorem]{Proposition}
\newtheorem{definition}[theorem]{Definition}
\newtheorem{corol}[theorem]{Corollary}

\newtheorem{remark}[theorem]{Remark}

\theoremstyle{definition}
\theoremstyle{remark}
\usepackage{xcolor}
\usepackage{hyperref}
\usepackage{lscape}

\begin{document}
	
\title[Generalized Exchange]{The discrete generalized exchange-driven system}

\author{P.K. Barik}
\address[P. K. Barik]{Instituto de Matem\'aticas,
Universidad de Granada, Rector L\'opez Argüet, 
S/N, 18001, Granada, Spain, 
Departamento de Matemática Aplicada,
Universidad de Granada,
Avenida de Fuentenueva S/N, 18071, Granada, Spain,
and Department of Mathematics, BITS-Pilani, Dubai Campus,
P.O. Box 345055, Dubai, United Arab Emirates
}  \email{barik@dubai.bits-pilani.ac.in}

\author{F.P. da Costa}
\address[F.P. da Costa]{Univ. Aberta, Dep. of Sciences and Technology,
  Rua da Escola Polit\'ecnica 141-7, P-1269-001 Lisboa, Portugal, and
  Univ. Lisboa, Instituto Superior T\'ecnico, Centre for Mathematical
  Analysis, Geometry and Dynamical Systems, Av. Rovisco Pais,
  P-1049-001 Lisboa, Portugal.}  \email{fcosta@uab.pt}

 \author{J.T. Pinto} 
\address[J.T. Pinto]{Univ. Lisboa, Instituto Superior T\'ecnico, Dep. of Mathematics
and Centre for Mathematical
Analysis, Geometry and Dynamical Systems, Av. Rovisco Pais,
P-1049-001 Lisboa, Portugal.}  \email{jpinto@tecnico.ulisboa.pt}

\author{R. Sasportes$^\dag$}
\thanks{$^\dag$ Our friend and co-worker Rafael Sasportes (1960--2024) died while 
this paper was on the final stages of preparation. We dedicate it to his memory.} 
\address[R. Sasportes]{Univ. Aberta, Dep. of Sciences and Technology,
Rua da Escola Polit\'ecnica 141-7, P-1269-001 Lisboa, Portugal, and
Univ. Lisboa, Instituto Superior T\'ecnico, Centre for Mathematical
Analysis, Geometry and Dynamical Systems, Av. Rovisco Pais,
P-1049-001 Lisboa, Portugal.}  \email{rafael.sasportes@uab.pt}

\thanks{Research partially supported by Funda\c{c}\~ao para a Ci\^encia e a Tecnologia (Portugal) 
through project CAMGSD UID/04459/2023.            }
\thanks{\emph{Corresponding author:} P. K. Barik}

\date{First version: September 7, 2024; \;Revised: July 17, 2025}

\subjclass{Primary 34A12; Secondary 34A34, 34A35, 92E20}

\keywords{exchange-driven cluster growth, aggregation kinetics, ordinary differential equations}

\begin{abstract}
We study a discrete model for generalized exchange-driven growth in which the particle exchanged 
between two clusters is not limited to be of size one. This set of models include as special cases the usual 
exchange-driven growth system and the coagulation-fragmentation system with binary fragmentation. Under 
reasonable general condition on the rate coefficients we establish the existence of admissible solutions, 
meaning solutions that are obtained as appropriate limit of solutions to a finite-dimensional truncation 
of the infinite-dimensional ODE. For these solutions we prove that, in the class of models we call isolated both the 
total number of particles and the total mass are conserved, whereas in those models we can non-isolated
only the mass is conserved. Additionally, under more restrictive growth conditions for the
rate equations we obtain uniqueness of solutions to the initial value problems.
\end{abstract}

\maketitle

\section{Introduction}\label{sec1}
The dynamics of growth processes of aggregates, or clusters,
is ubiquitous across all the natural world. For instance: in chemistry, 
polymerization processes \cite{vanDongen}; in astrophysics, the creation of stars and planets \cite{safronov}; 
in the physics of the atmosphere, the formation of clouds \cite{Pruppacher}; in solid state physics, 
deposition processes  \cite{mulheran}; 
in biology, the aggregation of red blood cells \cite{guy};in ecology, the grouping behavior of animals 
\cite{degond}; in economics, the merger of enterprises \cite{banakar}, and many more.  

\medskip

By assigning at each time a positive integer to each cluster size, representing their concentrations in the system, 
the evolution of the cluster sizes' distribution can be mathematically described by differential equations. 
One well-known mathematical 
model is the \textit{Smoluchowski coagulation equation}, see \cite{Smoluchowski:1917}, in which two clusters 
interact, leading to coalescence and the formation of a larger-sized cluster, and its generalization by inclusion 
of the possibility of fragmentation of the clusters (see, e.g.,
\cite{BLLvol1,BLL,dc15} and references therein). In the early 2000s several authors 
introduced a different growth model, known as the \textit{Exchange Driven Growth} model (EDG), 
arising, \emph{inter alia,} as the mean field limit for zero-range processes in non-equilibrium 
statistical physics, and also in modelling several social phenomena 
including migration, population dynamics and wealth exchange (see e.g. \cite{Ben:2003,IKR:1998,KL:2002,LR:2002}). 
The mechanism underlying this model involves the exchange of a single unit of mass (a \emph{monomer}, as it is 
usually called in the coagulation literature) from one cluster to another whenever they came into interaction.

\medskip 

The first mathematical study of the EDG equation was done by  Emre Esent\"urk in \cite{esen1} where the existence 
and uniqueness of solutions to the EDG equations is discussed. In addition, the gelation transition, instantaneous
gelation phenomenon and the existence of local solutions for different classes of interaction rates are also 
addressed. The large time behaviour of solutions was also investigated for different classes of reaction rates, 
see \cite{esen2,esenvel}.  In the recent work  \cite{schl1} Schlichting investigates the well-posedness of solutions to the 
EDG equation in which the conditions on initial data are relaxed when compared to those in \cite{esen1}.
Additionally, qualitative aspects such as the long-time behavior of the solution are also discussed. Furthermore, 
a study of self-similar solutions is carried out, specifically focusing on the product interaction rate, as detailed in
 \cite{Eichenberg:2021}. 
 
 \medskip
 
 In the present paper our primary goal is to propose a new, more general, discrete model for 
 exchange-driven growth of cluster and start its mathematical study. This new model 
 will be referred to as the \textit{discrete generalized exchange-driven growth} model (DGED, for short).
 As explained in the next session, the DGED model allows for the exchange not only of a monomer
 between two interacting clusters, but of the transfer of a bigger chunk of one cluster to the other.
 Thus, the DGED model includes as a special case the EDG model referred to above.
 In a certain sense, the usual coagulation-fragmentation equations (with binary fragmentation) can
 also be considered as a special case of the DGED equations.

\medskip

The outline of this article is as follows. In sections \ref{GEDGM} and \ref{MS}, we introduce the DGED 
model for two cases that, in section \ref{GEDGM}, are called \emph{isolated} and \emph{non-isolated}, and some relevant mathematical notions and results that will be needed later in this work.  In addition, all the main results of our paper are stated in section \ref{MS}. In section~\ref{sec2} 
we consider a truncated finite dimensional version of the DGED system. This is a finite dimensional ordinary differential equation system, for which the existence, uniqueness, and non-negativity of solutions to 
Cauchy problems are
automatically obtained from the standard theory. Here we prove some properties of its solutions
that will play an important role in the next section. In section \ref{sec3} we establish the existence of mild solutions to Cauchy problems for the  DGED system, i.e., continuous solutions of the integral version 
of the differential equations system. 
This is done by the use of Helly's theorem and a diagonal argument to establish
the existence of a function which is the limit of a (sub)sequence of solutions to truncated systems as
the dimension of the truncation grows to infinity, and then by proving that this function is indeed a (mild) solution of the DGED system using appropriate bounds (obtained in the previous two sections) and the bounded and monotone convergence theorems. Still in section~\ref{sec3} we prove that for the isolated case solutions
obtained in this way (which are called \emph{admissible}) conserve two important quantities, which
have the ``physical'' interpretation of being the total number of clusters, and the total mass
of the system, whereas in the non-isolated case only mass is conserved. In the short section \ref{secreg} we prove that under slightly more restrictive conditions
on the rate coefficients of the DGED system and assuming conservation of the total cluster mass (which
is true for admissible solutions in both the isolated and non-isolated cases),
then each component of the mild solution proved to exist in section~\ref{sec3} is indeed 
continuously differentiable.
In Section 7, we prove a uniqueness result under a more restrictive growth condition 
on the rate coefficients, applying standard techniques used in coagulation-type systems. 
We conclude the paper with a section of final remarks where we discuss some possible future
work related to the long-term behavior of solutions for this model.


\section{The generalized exchange-driven growth model}\label{GEDGM}
\begin{figure}[h]
	\psfrag{p}{$p$}
	\psfrag{q}{$q$}
	\psfrag{p-k}{$p-k$}
	\psfrag{q+k}{$q+k$}
         \includegraphics[scale=0.70]{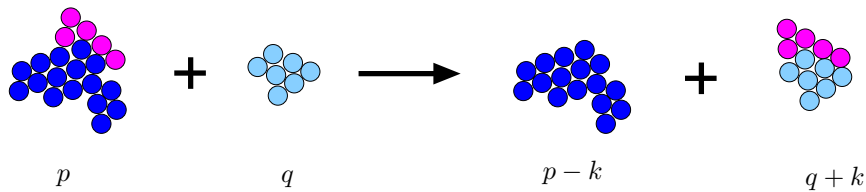}
	\caption{Reaction scheme of the DGED model considered in this paper. 
A chunk of size $k$ in a particle of size $p$ is transferred to a particle of size $q$ to produce a
particle of size $p-k$ and another of size $q+k$. The rate coefficient for this reaction is $a(p,q; k)\geq 0.$}
\label{fig1}
\end{figure}

Consider a population of particles with sizes described by $p\in\Nb = \{0, 1, 2, \ldots\}.$ For each size $p$ and each time $t$, 
let $c_p(t)$ be the concentration of particles of size $p$ (or $p$-clusters, for short)
at time $t$. Then, we assume that a chunk of size $k\leq p$ can be detached from this particle and attached to another of size $q$, schematically represented by
\begin{equation}
\langle p\rangle + \langle q\rangle \rightarrow \langle p-k\rangle + \langle q+k\rangle.\label{eqscheme}
\end{equation}
Figure~\ref{fig1} gives a pictorial illustration of this process.

%
%
%
%
%
%

\medskip

Assuming the validity of the mass action law of chemical kinetics, the rate of production of  $(p-k)$-clusters and $(q+k)$-clusters due to the
reaction between $p$-clusters and $q$-clusters, according to the mechanism displayed in Figure~\ref{fig1}, is equal to $a(p,q; k)c_{p}(t)c_q(t),$ where the rate coefficient is $a(p,q; k).$

\medskip

Hereafter $a(p, q; k)$ will denote the rate of the reaction in which a $k$-cluster is detached from a $p$-cluster and attaches itself to a $q$-cluster. 
Clearly, we must have $p, q, k\in \Nb$ with $p\geq 1$, and $1 \leq k \leq p$ because the case $k=0$ is the absence of reaction and $k>p$ does not make physical
sense as we cannot detach from a given cluster a part bigger than the entire cluster.
However, we will allow $k=p$ in \eqref{eqscheme}, which means that the entire $p$-cluster is attached to the $q$-cluster in the reaction
\[
\langle p\rangle + \langle q\rangle \rightarrow \langle 0\rangle + \langle q+p\rangle.
\]
This is very much like the usual coagulation reaction (see, e.g., \cite{dc15}) but for the consideration of the ``void'', or
``empty'',  cluster $\langle 0\rangle$.

\medskip

\begin{figure}[h]
	\psfrag{p}{$p$}
	\psfrag{q}{$q$}
	\psfrag{r}{$k$}
	\psfrag{r=p}{$k=p$}
	\psfrag{q=0}{$q = 0$}
	\includegraphics[scale=0.35]{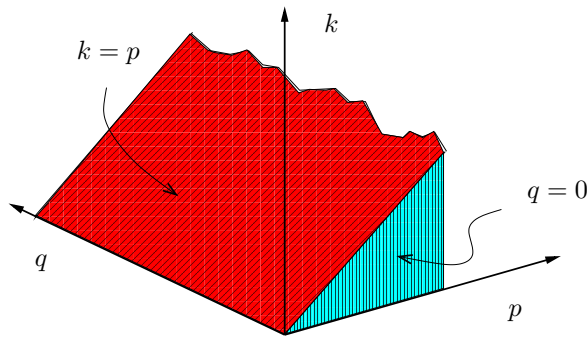}
	\caption{Domain of the rate coefficients $a(p, q; k)$: region of $\Nb^{3}$ bounded by the planes
$q=0$ (corresponding to binary fragmentation) and $k=p$ (corresponding to coagulation).}\label{fig2}
\end{figure}

In a similar way, if we consider $q=0$ in the reaction scheme \eqref{eqscheme} it becomes
\[
\langle p\rangle + \langle 0\rangle \rightarrow \langle p-k\rangle + \langle k\rangle,
\]
and this is a kind of fragmentation of the non-void cluster. 
Figure~\ref{fig2}  illustrates the domain of $a(p, q; k)$ and the regions just described.

%
%
%
%
%
%

\medskip

Allowing the kind of reactions
where the void clusters are destroyed by reaction with other clusters, 
means that one is considering the class of so called active exchange
driven models (in the classification of Esenturk and  Connaughton \cite{esen2}). This is what
will be done in this paper.    

\medskip

The way the dynamics of the void clusters is considered provides a further classification of the 
models under consideration: two natural choices correspond to open, or \emph{non-isolated}, and to \emph{isolated} thermodynamic
systems. In the first case we can assume that the cluster system is embedded in an infinite size bath of
void particles such that the concentration of $\langle 0\rangle$ particles remains constant in time, 
i.e., $\dot{c}_0(t)=0, \forall t$; in the second case one
assumes that the dynamics of the cluster system, including the void particles, is such that the total 
number of clusters 
is conserved (at least at a formal level)
and so there is the need to have a dynamic equation describing the evolution of the concentration 
of clusters $\langle 0\rangle$ compatible with this assumption. 
This same distinction can be see in \cite{penrose89} for the
role of the monomers ($\langle 1\rangle$ particles) in the context of the 
dynamic of  cluster equations of Becker-D\"oring type.

\medskip

Returning to the above reaction schemes involving the void cluster, we note that its presence can
have two distinct effects according to whether the system is isolated or not. In the non-isolated case, where the concentration of $\langle 0\rangle$ is constant,
the reactions are exactly the coagulation and the binary fragmentation 
reactions considered in coagulation-fragmentation studies.
Thus, for that case, the rate coefficients $a(p,q;p)$ and $a(p,0,k)$ correspond to reactions of 
coagulation and binary fragmentation, respectively, and all the other
cases (with $0<k< p$) correspond to genuine exchange of pieces of the $p$-cluster to the $q$-cluster.
In the case of the isolated model, where the concentration of $\langle 0\rangle$ is not necessarily constant,
then the mass action law for the reaction corresponding to the destruction of void clusters 
produces a nonlinear differential equation and the parallel with the usual
fragmentation equations is less straightforward.

\medskip

The DGED system describes, for each $i\in\Nb,$ the time evolution  of the concentration of the $i$-clusters and is obtained by keeping track of the various ways an $i$-cluster can be formed or destroyed by the kinetics of the type \eqref{eqscheme}. As in previous studies of exchange-driven kinetics, such as \cite{esen1,esenvel,schl1}, we also consider clusters of size zero, the void, or empty, clusters. As we shall see this is rather convenient, enabling, in the isolated model, the existence of an additional conservation law, and, in the non-isolated model, the
recovery of the usual discrete coagulation-fragmentation model as a particular case.

\medskip

In order to introduce the DGED system
we have to consider the various ways that  lead to the creation and the 
destruction of an $i$-cluster, which we do next, starting with the case of the isolated system.

\medskip

In scheme \eqref{eqscheme}, one $i$-cluster is created if (i) $p-k=i$, or if (ii) $q+k=i.$ Likewise, one $i$-cluster is destroyed if (iii) $q=i$, or if (iv) $p=i.$ In the case of creation (i), for every $k\in\Nb_+:=\Nb\setminus\{0\},$ we have to consider all the indices combinations for which $p=i+k$ and $q=j,$ for any $j\in\Nb.$
Therefore, in this case, the formation of an $i$-cluster proceeds through the scheme \eqref{eqscheme} illustrated in Figure~\ref{fig1} and, 
by the mass action law, its contribution to the global rate of change of $c_i(t)$ is given by
\begin{equation}
Q_{1, i}(c(t)) := \sum_{k=1}^\infty \sum_{j=0}^\infty a(i+k, j; k)c_{i+k}(t)c_{j}(t).		\label{Q1i}
\end{equation}
For the destruction case (iii), we have to consider all possible indices combinations with $q=i,$ that is, $j:=p\geq k\geq 1$. The corresponding contribution to the global rate of change of $c_i(t)$ is thus
\begin{equation}
Q_{2, i}(c(t)) := - \sum_{k=1}^\infty \sum_{j=k}^\infty a(j, i; k)c_{j}(t)c_{i}(t).		\label{Q2i}
\end{equation}
For the creation case (ii), we have to consider all possible indices combinations for which $q=i-k,$ that is, for $1\leq k\leq i,$ we have $p\geq k,$ which, by calling $j:=p$ leads to the contribution
\begin{equation}
Q_{3, i}(c(t)) := \sum_{k=1}^i \sum_{j=k}^\infty a(j, i-k; k)c_{j}(t)c_{i-k}(t).		\label{Q3i}
\end{equation}
Finally, for the destruction case (iv), we have to consider the indices combinations for which
$k\leq p=i,$ $q\in\Nb.$ By calling $j:=q,$ the corresponding contribution to the rate equation is
\begin{equation}
Q_{4, i}(c(t)) := -  \sum_{k=1}^i \sum_{j=0}^\infty a(i, j; k)c_{j}(t)c_{i}(t).		\label{Q4i}
\end{equation}

\medskip

Clearly the processes leading to $Q_{3, i}$ and $Q_{4, i}$ cannot be in operation when $i=0$ because they require the consideration of clusters with 
size smaller than $i$. Hence we define $Q_{3, 0} = Q_{4, 0}=0.$

\medskip

If the system is not isolated but instead the number density of $0$-clusters is kept constant at some fixed
value $c_0(t) = c_{00}$ the equation for the evolution of $c_0(t)$
becomes simply
\begin{equation}
 \dot{c}_0 = 0.    \label{DGEDnonisolated}
\end{equation}
Hence, for this case $Q_{j,0}=0,$ for all $j=1, 2, 3, 4.$

\medskip
Thus, in all cases, the Discrete Generalized Exchange Driven system (DGED) can be written as 
\begin{equation}
\dot{c}_i = \sum_{j=1}^{4} Q_{j, i}(c), \qquad i\in\Nb.    \label{DGED}
\end{equation}

We are particularly interested in the study of Cauchy problems for \eqref{DGED} considering the initial condition 
\begin{equation}
c_i(0) = c_{0 i}\geq 0,\quad\text{for $i\in \Nb.$} \label{init}
\end{equation}

\medskip

We remark that all terms $Q_{j,i}$ contain terms of the type $a(p, 0; p) c_{p}(t)c_{0}(t),$ which
corresponds to reactions
\[
\langle p\rangle + \langle 0\rangle \rightarrow \langle 0\rangle + \langle p\rangle,
\]
that, in fact, do not correspond to any change in the cluster size distribution. 
They are included in each of the $Q_{j,i}$ for notational convenience only and, in fact,
they cancel each other in \eqref{DGED} 
(the $j=0$ terms in $Q_{1,0}$ cancels the $j=k$ terms in $Q_{2,0}$, and the $k=i=j$ term in $Q_{3,i}$ cancels the term with
$k=i$ and $j=0$ in $Q_{4,i}.$) 
However, when taken in isolation each $Q_{j,i}$ contains these spurious, unphysical, 
contributions and since in the definition of 
solution one requires that each  $Q_{j,i}(c(\cdot))$ be integrable in $(0, T)$ (see Definition~\ref{Def Sol}) 
we need to make some assumption on the coefficients $a(p,0; p)$.
Given the discussion above we define them as
\begin{align}
& a(p, 0; p) = 0,\; \forall p \in \Nb_+.\label{nullcoef} 
\end{align}

\noindent
We will also assume that the following symmetry relations always hold:
\begin{align}
& a(k, j; k) = a(j, k; j),\; \forall j, k \in \Nb_+, \label{symcoag} \\
& a(i, 0; k) = a(i, 0; i-k),\; \forall k\in \{1, \ldots, i-1\}, i\in\Nb_+. \label{symfrag}
\end{align}
Condition \eqref{symcoag} is due to the fact that, as pointed out above, each rate coefficient $a(i, j;k)$ with $i=k$ corresponds to the rate of a coagulation reaction of a $k$-cluster with a $j$-cluster
(to produce a $(j+k)$-cluster) and thus this is the usual symmetry assumption due to the fact that a coagulation
reaction of a $k$ with a $j$-cluster is the same as the reaction of a $j$ with a $k$-cluster. Similarly, \eqref{symfrag} 
correspond to the symmetry of the rates of the reaction
$\langle i\rangle + \langle 0\rangle \rightarrow \langle i-k\rangle + \langle k\rangle$, 
when interpreted as a fragmentation of the $i$-cluster produced by the shedding off a $k$-cluster, or by the shedding off a $(i-k)$-cluster.

\medskip

Observe that if the rate coefficients satisfy $a(i,j;k)=0$ if $k\not= 1,$ we obtain the exchange-driven growth system introduced in \cite{Ben:2003} and  studied  mathematically in \cite{Eichenberg:2021,esen1,esen2,esenvel,schl1}, 
with rate coefficients $K(i,j).$ In this case, 
\begin{equation}\label{EDG1}
\forall i,k\in\Nb_+, \; j\in \Nb,\quad  a(i,j;k)=K(i,j)\delta_{k,1},
\end{equation}
where $\delta_{k,1}$  is the Kronecker $\delta$-symbol equal to $1$ if $k=1$ and zero otherwise,.

\medskip

As was already pointed out above, it is interesting to observe that the DGED system formally 
reduces to the standard discrete coagulation-fragmentation equations in the non-isolated case, provided the 
rate coefficients satisfy the conditions
\begin{align}
& a(k, j; k)\geq 0,  \label{coag}\\
& a(i, 0; k)\geq 0,  \label{frag}\\
& \text{and all the other cases are equal to zero,} \nonumber
\end{align}
In fact, observe that with these conditions we get
\begin{align}
Q_{1,i}(c) & 
 \bigeq{\eqref{frag}} \sum_{k=1}^\infty  a(i+k, 0; k)c_{i+k}c_{0} 
 \label{q1}\\
Q_{2, i}(c) & 
\bigeq{\eqref{coag}} - \sum_{k=1}^\infty a(k, i; k)c_{k}c_{i} \label{q2}\\
Q_{3, i}(c) & 
\bbigeq{\eqref{coag}, \eqref{frag}}  \sum_{j=i}^\infty a(j,0; i)c_jc_0 + \sum_{k=1}^i a(k, i-k; k)c_kc_{i-k} \label{q3} \\
Q_{4, i}(c) &
\bbigeq{\eqref{coag}, \eqref{frag}} - \sum_{k=1}^i a(i,0; k)c_ic_0 - \sum_{j=0}^\infty a(i, j; i)c_ic_{j} \label{q4}
\end{align}
Now, due to \eqref{nullcoef}, the terms $j=i$ in the first sum in \eqref{q3}, $k=i$ in the second sum in \eqref{q3}
and in the first sum in \eqref{q4}, and $j=0$ in the second sum in \eqref{q4} are all equal to zero, and using 
the symmetry conditions \eqref{symcoag} and \eqref{symfrag} 
\begin{align}
\sum_{j=1}^\infty a(i, j; i)c_ic_j & \bigeq{\eqref{symcoag}} \sum_{j=1}^\infty a(j, i; j)c_ic_j \label{qq4}\\
\sum_{j=i+1}^\infty a(j,0; i)c_jc_0 & \bigeq{\eqref{symfrag}} \sum_{j=i+1}^\infty a(j,0; j-i)c_jc_0 = \sum_{k = 1}^\infty a(i+k,0; k)c_{i+k}c_0,\label{qq3}
\end{align}
then, putting together \eqref{q1}-\eqref{qq4} we obtain, for $i>1$, 
\begin{align}
\sum_{j=1}^{4} Q_{j, i}(c) 
& = \sum_{k=1}^{i-1}a(k, i-k; k)c_kc_{i-k} - \sum_{k=1}^{i-1}a(i, 0;k)c_ic_0 \, + \nonumber \\
& \;\;\; + 2 \sum_{j=1}^\infty a(i+j, 0; j)c_{i+j}c_0 - 2\sum_{j=1}^\infty a(i,j; i)c_ic_j  \nonumber 
\end{align}
and so we can write the DGED system  as
\begin{align}
\dot{c}_i & = \frac{1}{2}\sum_{k=1}^{i-1}W_{k, i-k}(c) - \sum_{k=1}^\infty W_{i,k}(c),
\qquad i\in\Nb_+ \label{usualcf}
\end{align}
where $W_{i,j}(c):= 2a(i, j; i)c_i c_j - 2a(i+j, 0; j)c_{i+j}c_0$. 
Clearly, if our system is immersed in an infinite particle bath of $0$-cluster particles (which means that the concentration of $0$-particles is kept
constant) and if we define $a_{i,j}:=2 a(i, j; i)$ and $b_{i,j}:= 2a(i+j, 0; j)c_0$, then \eqref{usualcf} becomes the usual discrete
coagulation-fragmentation equation \cite{dc15}. Therefore, in this case we have
\begin{align}
\forall i,j,k\in\Nb_+,\;k\leq i,\quad a(i,j;k)=\frac{1}{2}a_{i,j}\delta_{k,i},\label{EDGCF1}\\
\forall i,k\in\Nb_+,\; k\leq i, \quad a(i,0;k)c_0=\frac{1}{2}b_{i-k,k}.\label{EDGCF2}
\end{align}



\section{Mathematical setting and main results}\label{MS}
In this section we introduce the notion of solution to the initial value problem \eqref{DGED}, \eqref{init}
 that we consider in this work as well as other concepts that we will be useful in the following. Also,  we introduce here the hypothesis we will consider for the rate coefficients $a(i,j;k).$ 

\medskip

We start by introducing the concentration moments. For $r \geq 0$, we denote by $\mathcal{P}_r(t)$ the $r$-th moment of the solution $c(t)$ at a given time $t$:
\begin{align}
\mathcal{P}_r(t) := \sum_{i=0}^{\infty} i^r c_i(t).   \label{moments}
\end{align}
In the cases of $r=0$ and $r=1$, the moments have a natural physical interpretation: $ \mathcal{P}_0(t)$ represents the total number of particles in the system, while  $\mathcal{P}_1(t)$ is the system total mass. 
This is the reason that, for these type of models, it is anticipated that both these moments 
are invariant under time evolution.
\medskip

The above-mentioned physical interpretation suggests that, analogously to previous works on 
coagulation-type systems, it is natural to work in the Banach space
\begin{equation}
X_{0, 1} := \Bigl\{c:\Nb\to\Rb \,\bigl|\; \|c\| := \sum_{i=0}^\infty (1+i)|c_i|  < \infty\Bigr\}.		\label{X01}
\end{equation}
Note that we can write $\|c\| = \|c\|_{\ell_1} + \|(ic_i)\|_{\ell_1}$ where  $\|u\|_{\ell_1}$ is the usual $\ell_1$ norm of the sequence $u=(u_i):\Nb\to\Rb$.
\medskip

Also due to the physical interpretation of $c=(c_i)$ as a sequence of concentrations, we will be exclusively interested in
solutions in the non-negative cone of $X_{0, 1}$,
\begin{equation}
X_{0, 1}^+ := X_{0, 1} \cap \{c=(c_i)\,|\; c_i\geq 0\}.			\label{X01+}
\end{equation}
 
 \medskip
 
We are now ready to state the following definition of solution.
 \begin{definition}\label{Def Sol}
Let $T\in (0, +\infty]$ and let $c_0 = (c_{0 i}) \in X_{0,1}^+$ be a sequence of nonnegative real numbers. A
(mild) solution to \eqref{DGED}, \eqref{init}
on $[0, T)$ is a sequence of nonnegative continuous functions
$c=(c_i): [0, T)\to X_{0,1}^+$, such that, for each $i\in\Nb$ and $t\in (0, T)$, the following holds:
\begin{itemize}
\item[(i)] $c_i\in C^0([0, T)),$
\item[(ii)] $Q_{j, i}(c(\cdot)) \in L^1(0,t),\;\, j\in \{1, 2, 3, 4\}$,
\item[(iii)] and
\[
c_i(t) = c_{0 i} + \int_0^t\sum_{j=1}^{4} Q_{j, i}(c(s)) \,ds.
\]
\end{itemize}
\end{definition}

\subsection{Rate Coefficients: Growth Assumptions and Examples}
As stated previously, we assume 
conditions \eqref{nullcoef}, \eqref{symcoag} and \eqref{symfrag} throughout the paper.
With respect to the bounds on the rate coefficients, we assume the existence of positive constants $C$ and $\mathcal{Q}$, which as a matter of convenience we take $C,\mathcal{Q}\geq 1$, as well as non-negative numbers $q_{i,k},$ with
$1\leq k\leq i,$ satisfying
\begin{equation}\label{qikbound}
\forall i\in\Nb_+,\quad \sum_{k=1}^i k(i-k+1)q_{i,k}\leq \mathcal{Q}i,
\end{equation}
such that, for all integers $i,k$ such that $1\leq k\leq i,$ and all $j\in\Nb_+,$
\begin{equation}\label{bound}
 a(i,j;k)\leq C(i-k+1)(j+k)q_{i,k}.
\end{equation}
We can see $q_{i,k}$ as measuring the ease with which a cluster of size $k$ can be detached from a
cluster of size $i$: the smaller the value of $q_{i,k}$ is, the smaller is the rate coefficient $a(i,j; k)$, and so the
harder it is for a $k$-cluster to take part in the exchange reaction $\langle i\rangle + \langle j\rangle \to  \langle i-k\rangle + \langle j+k\rangle.$

\medskip

One important bound for our subsequent arguments comes from the fact that, for every pair of integers $i,k$ such that $1\leq k\leq  i,$ we have
$$
i\sum_{k=1}^i q_{i,k} = \min_{1\leq k\leq i}\{k(i-k+1)\}\sum_{k=1}^i q_{i,k}\leq\sum_{k=1}^i k(i-k+1)q_{i,k}\leq \mathcal{Q}i,
$$
thus resulting in ${\displaystyle\sum_{k=1}^i q_{i,k} \leq \mathcal{Q}}$, and, in particular,
\begin{equation}\label{qikest}
q_{i,k}\leq \mathcal{Q}.
\end{equation}

\medskip

We illustrate next by exhibiting some examples that \eqref{qikbound} and \eqref{bound}
are fulfilled by some well known cases in the literature.

\medskip

\noindent
{\bf Example 1.\/} Consider
 the exchange driven growth system with kernel satisfying $K(i,j)\leq C_0ij$, for $i,j\in\Nb_+,$ by \eqref{EDG1}, we have, for all $i,j\in\Nb_+$
$$
a(i,j;k)=K(i,j)\delta_{k,1}\leq C_0ij\delta_{k,1}\leq C_0i(j+1)\delta_{k,1} 
$$
which is \eqref{bound} with $C=C_0$ and $q_{i,k}=\delta_{k,1},$ observing that \eqref{qikbound} is also satisfied since,
$$
\sum_{k=1}^i k(i-k+1)q_{i,k}=\sum_{k=1}^i k(i-k+1)\delta_{k,1}=i,
$$
so that we can take $\mathcal{Q}=1$ in this case.

\medskip

\noindent
{\bf Example 2.\/} 
In fact, Example 1 is a particular case of the more general case where there is an upper bound, $\bar{k}\geq 1,$ for the number of particles exchanged between two reacting clusters of any sizes, that is,
$$
a(i,j;k)=0,\quad\text{if}\quad k>\bar{k}.
$$
We now show that, for this case, the existence of a constant $\bar{C}>0$ such that, for $i,j,k\in\Nb_+,$
\begin{equation}\label{kbarcond}
a(i,j;k)\leq \bar{C}ij,\qquad 1\leq k\leq \min(i,\bar{k})
\end{equation}
is equivalent to the existence of $C\geq 1$ such that \eqref{bound} is satisfied with $q_{i,k}=\mathds{1}_{\{1,\dots,\bar{k}\}}(k),$ for which \eqref{qikbound} is true.  In first place, it is easy to see that this choice of $q_{i,k}$ verifies \eqref{qikbound}:
$$
\sum_{k=1}^i k(i-k+1)\mathds{1}_{\{1,\dots,\bar{k}\}}(k)=\sum_{k=1}^{\min(i,\bar{k})} k(i-k+1)\leq \bar{k}^2i.
$$
Furthermore, on one hand, for any $i,j,k\in\Nb_+,$ such that $1\leq k\leq \min(i,\bar{k}),$ we have
$$
(i-k+1)(j+k)\leq ij\left(1+\frac{k}{j}\right)\leq (1+\bar{k})ij,
$$
and on the other, we also have, in case $i\geq \bar{k},$
$$
(i-k+1)(j+k)\geq ij\left(1-\frac{k-1}{i}\right)\geq ij\left(1-\frac{\bar{k}-1}{\bar{k}}\right)=\frac{1}{\bar{k}}ij,
$$
and in case $i<\bar{k},$
$$
(i-k+1)(j+k)\geq j\geq\frac{1}{\bar{k}}ij.
$$
Therefore, if condition \eqref{bound} holds with the above choice of $q_{i,k}$ which also verifies \eqref{qikbound}, then condition \eqref{kbarcond} is verified with $\bar{C}=(1+\bar{k})C.$ Reciprocally, if the rate coefficients verifies \eqref{kbarcond} then conditions \eqref{qikbound},\eqref{bound} hold with the above choice of $q_{i,k}$ and $C=\bar{k}\bar{C}.$

\medskip

\noindent
{\bf Example 3.\/} 
An example of a strictly positive
rate coefficient $a(i,j;k)$ without 
a fixed upper bound $\bar{k}$ on the size of clusters 
being exchanged is the constant kernel $a(i,j;k)=1$, for
which we have
\[
1 = \frac{(i-k+1)(j+k+1)}{(i-k+1)(j+k+1)} \leq 
2\,\frac{(i-k+1)(j+k)}{(i-k+1)k} =:2\,(i-k+1)(j+k)q_{i,k}.
\]
Clearly, the symmetry conditions \eqref{symcoag} and
\eqref{symfrag}, as well as the bound \eqref{qikbound}, 
are trivially satisfied.

\medskip

\noindent
{\bf Example 4.\/} 
A less trivial example of the sort considered in Example 3 but such that the
rate coefficients are unbounded is
\[
a(i,j;k) = \frac{(i-k+1)(j+k+1)}{1 + (i-k)k},
\]
for which checking the symmetry and growth conditions
are also easily done.

\medskip

\noindent
{\bf Example 5.\/} 
Finally, in this example let us consider  the coagulation-fragmentation case 
with upper bound on the coagulation rate coefficients of the type $a_{i,j}\leq C_1(i+j),$ by \eqref{EDGCF1}. 
For this case we have
$$
a(i,j;k)=\frac{1}{2}a_{i,j}\delta_{i,k}\leq\frac{C_1}{2}(i+j)\delta_{i,k},
$$
which is again \eqref{bound}, this time taking $C=\max(1,C_1/2)$ and $q_{i,k}=\delta_{i,k},$ remarking also that
$$
\sum_{k=1}^i k(i-k+1)q_{i,k}=\sum_{k=1}^i k(i-k+1)\delta_{i,k}=i
$$
so that we can again take $\mathcal{Q}=1.$ Remark that no bounds are imposed on the coefficients $a(i,0;j),$ which allows us to consider any choice of fragmentation coefficients $b_{i-j,j}.$ 

\medskip

Therefore, from these examples (in particular from examples 1 and 5) we conclude that the cases of the exchange driven system growth with kernel of multiplicative bound type and the coagulation-fragmentation system with coagulation kernel with linear additive bound 
are both particular cases of the discrete generalized exchange driven system with bounds of the type \eqref{qikbound}-\eqref{bound}.

\medskip

\subsection{Main results}
In this subsection, we state the main results proved in this paper
for the initial value problem \eqref{DGED}–\eqref{init}, considering both the isolated and non-isolated cases.

Let us first state the existence result for the isolated case.
\begin{theorem}\label{ExistThm1}
Let $c_0 =  (c_{0 i})_{i \ge 0}  \in X_{0, 1}^+ $. Assume the rate coefficients satisfy the conditions \eqref{qikbound}--\eqref{bound}. Then the initial value problem \eqref{DGED}, \eqref{init} for the isolated case has at least one global solution $c$ defined on the interval $ [0, \infty)$. 
\end{theorem}

The next theorem establishes mass and particle number conservation in the isolated system.
\begin{theorem}\label{TPTheorem}
Let $c_0\in X_{0,1}^+$, and let $c= (c_i )_{i \geq 0}$ be an admissible solution to 
\eqref{DGED}, \eqref{init} in the isolated case. Then, for all $t\geq 0$,
\begin{align}\label{P0conserve}
\mathcal{P}_0(t) := \sum_{i=0}^{\infty} c_i(t) = \sum_{i=0}^{\infty} c_{0 i} =: \mathcal{P}_0(0).
\end{align} 
and
\begin{align}\label{P1conserve}
\mathcal{P}_1(t) := \sum_{i=0}^{\infty} i c_i(t) = \sum_{i=0}^{\infty} i c_{0 i} =: \mathcal{P}_1(0).
\end{align} 
\end{theorem}

Similarly to Theorem~\ref{ExistThm1}, we now state the existence of global solutions for the non-isolated case.
\begin{theorem}\label{ExistThm2}
Let $c_0 =  (c_{0 i})_{i \ge 0}  \in X_{0, 1}^+ $. Assume
the rate coefficients satisfy the conditions \eqref{qikbound}--\eqref{bound}. Then the initial value problem \eqref{DGED}, \eqref{init} for the non-isolated case has at least one global solution $c$ defined on the interval $ [0, \infty)$.
\end{theorem}

The following result establishes the mass conservation property in the non-isolated case:
\begin{theorem}\label{TPTheorem2}
Let $c_0\in X_{0,1}^+$, and let $c= (c_i )_{i \geq 0}$ be an admissible solution to 
\eqref{DGED}, \eqref{init} in the non-isolated case. Then, for all $t\geq 0$,
$\mathcal{P}_1(t) = \mathcal{P}_1(0).$
\end{theorem}

We next state the regularity of the solution in classical sense: 
\begin{theorem}[Regularity of solutions]\label{Regthm}
Suppose that the rate coefficients satisfy \eqref{qikbound}-\eqref{bound}, for all $j\in\Nb.$ If $c=(c_i)$ is a
mild solution of \eqref{DGED}-\eqref{init} on $[0,T)$, 
in either the isolated or the non-isolated cases, and if 
$\sum_{i=1}^\infty ic_i(\cdot)$ is constant in $[0,T),$ 
then, for each $i\in \Nb,$  $c_i\in C^1([0,T),\mathbb{R}^+_0).$
\end{theorem}

We finally state a partial uniqueness result:
\begin{theorem}\label{Uniqthm}
Suppose that the rate coefficients satisfy \eqref{bound} and \eqref{boundalpha}, with $\alpha \in [0, \frac{1}{2})$. 
Let $T \in (0, +\infty)$, and $c_0 \in X_{0,1}^+$. Then, the initial value problem \eqref{DGED}--\eqref{init} has
a unique solution on $[0, T]$.
\end{theorem}

\subsection{Preliminary results}
Prior to establish the existence result in section~\ref{sec3}, it is imperative to revisit certain notations
and to get some inequalities that will be needed later. Let us first recall the space $\mathcal{E}$ from \cite{PL:2002}:
\begin{quote}
	$\mathcal{E}$ is the set of non-negative 
	and convex functions $\sigma \in \mathcal{C}^1([0, \infty)) \cap W^{2, \infty}_{loc}(0, \infty) $ 
	with $\sigma(0)=0$, $\sigma'(0) \geq 0$ and $\sigma'$ is a concave function.
	Assume furthermore 
	that these functions satisfy the following condition:
	\begin{equation}\label{Convex}
		\lim_{r \to \infty } \sigma'(r) =   \lim_{r \to \infty } \frac{\sigma(r)}{r} = \infty.
	\end{equation}
\end{quote}	
Again from \cite{PL:2002}, define $\mathcal{E}_1$ as follows:
\begin{quote}
		$\mathcal{E}_1$ is the  set of all non-negative and convex functions 
 $\sigma \in \mathcal{C}^2([0, \infty))$ with $\sigma(0)=0$, $\sigma'(0) = 0$, and $\sigma'$ is 
 a convex function satisfying the so-called $\Delta_2-$condition, namely, there is a constant $ A_{\sigma} \geq 0,$ 
 such that 
 \begin{align}\label{Convex2}
  \sigma'(2x) \leq A_{\sigma} \sigma'(x), \ \ \ x \in [0, \infty). 
 \end{align}
\end{quote}

As an illustration of these concepts consider the following example: if $\sigma(x) = x^{p}$, then $\sigma \in \mathcal{E}$ when
$p \in (1, 2]$, and $\sigma \in \mathcal{E}_1$ if $p \geq 2$.

  \medskip
  
Instrumental to the proofs of our subsequent results is the following inequality that we recall from \cite[Lemma  3.2.]{PL:2002}, which is valid for all $\sigma\in \mathcal{E}\cup\mathcal{E}_1:$
 \begin{align}\label{Convex1}
(i+j)(\sigma(i+j)  - \sigma(i) - \sigma(j) )\leq m_{\sigma}  (i \sigma(j) + j \sigma(i)) , \quad \forall i \geq 0, j \geq 1.
 \end{align} 
 Here, $m_\sigma\geq 0$ is a constant that, in case $\sigma\in\mathcal{E},$ can be taken as $m_\sigma=2.$ 
 We have, as a consequence,
 \begin{lemma}\label{sigmalemma}
 Assume $\sigma\in\mathcal{E}_1\cup\mathcal{E}.$ Then, for all $i,j,k\in\Nb_+$ such that $k\leq i,$
 we have
 $$
 (j+k)(\sigma(i-k)+\sigma(j+k)-\sigma(i)-\sigma(j))\leq m_\sigma ( j \sigma(k) + k \sigma(j) ).
 $$
 \end{lemma}
 \begin{proof}
 First, write
 \begin{multline}\label{sijk}
 \sigma(i-k)+\sigma(j+k)-\sigma(i)-\sigma(j)\\= \Big(\sigma(j+k)-\sigma(j)-\sigma(k)\Big)-\Big(\sigma(i)-\sigma(i-k)-\sigma(k)\Big).
 \end{multline}
For each $a\geq 0$ define $\phi_a(x):=\sigma(x)-\sigma(a)-\sigma(x-a).$ Then, for any $x\geq a,$ by the convexity of $\sigma$ we have,
$$
\phi_a'(x)=\sigma'(x)-\sigma'(x-a)\geq 0.
$$
But then, since $\phi_a(a)=0,$ we conclude that, $\phi_a(x)\geq 0.$ Therefore,
$$
\sigma(i)-\sigma(i-k)-\sigma(k)=\phi_{i-k}(i)\geq 0.
$$
By using this in \eqref{sijk}, together with \eqref{Convex1} the proof is completed.
 \end{proof}

We also have the following important property:

 \begin{lemma}\label{sigmalemma2}
 Let $\sigma\in\mathcal{E}.$ Then, there are constants $\eta>0$ and $M_0>0$ such that, for any integers $p\geq M_0,$ and $k\in [1,p-1],$
 $$
 \sigma(p)-\sigma(p-k)-\sigma(k)\geq \eta\frac{\sigma(p-1)}{p-1}.
 $$ 
 \end{lemma}
 \begin{proof}
 Define $$\phi_p(k):=\sigma(p-k)+\sigma(k),\quad k\in[1,p].$$
 Then, $\phi_p(k)=\phi_p(p-k)$ and also, if $1\leq k\leq p/2,$ then $p/2\leq p-k\leq p-1.$ Hence, by the convexity of $\sigma,$ we have that $\sigma'(p-k)\geq \sigma'(k),$ and therefore, 
$$
\phi_p'(k)=-\sigma'(p-k)+\sigma'(k)\leq 0,\qquad \phi_p'(p-k)\geq 0.
$$
This implies that $\max_{k\in[1,p-1]}\phi_p(k)=\phi_p(1)=\phi_p(p-1).$ Thus,
\begin{equation}\label{minsigmadif}
\min_{k\in[1,p-1]}\big(\sigma(p)-\sigma(p-k)-\sigma(k)\big)=\sigma(p)-\sigma(p-1)-\sigma(1).
\end{equation} 
By the Lagrange theorem and the convexity of $\sigma$ we can infer
$$
\sigma(p)-\sigma(p-1)\geq \sigma'(p-1).
$$  
On the other hand, again by the Lagrange theorem and convexity of $\sigma,$
$$
\sigma(p-1)=\sigma(p-1)-\sigma(0)\leq \sigma'(p-1)(p-1),
$$
and therefore,
$$
\sigma(p)-\sigma(p-1)-\sigma(1)\geq \frac{\sigma(p-1)}{p-1}-\sigma(1)
=\frac{\sigma(p-1)}{p-1}\Big(1-\sigma(1)\frac{p-1}{\sigma(p-1)}\Big).
$$
Fix $\eta\in(0,1).$ Since $\sigma\in\mathcal{E},$ there is $M_0>0$ such that, 
$$
p\geq M_0\quad\Longrightarrow\quad \sigma(1)\frac{p-1}{\sigma(p-1)}\leq 1-\eta,
$$
so that,
$$
\sigma(p)-\sigma(p-1)-\sigma(1)\geq \eta\frac{\sigma(p-1)}{p-1}.
$$
By using this inequality in \eqref{minsigmadif} the lemma is proved.
\end{proof} 
 
%

%
%
\section{The truncated system}\label{sec2}

We will consider $N$-truncated systems that are obtained from \eqref{DGED} assuming that no clusters with sizes bigger than $N$ exist initially neither can they be formed by time evolution. This corresponds to a modification of the rate coefficients by making them zero
whenever the reaction in question involves clusters larger than $N$. An equivalent and maybe a more transparent way of 
reflecting this is to appropriately modify the sums in $Q_{j, i}$ and consider the finite-dimensional system with $i\in \{0, 1, \ldots, N\}.$ It is not hard to see that, in the case of isolated systems, the following system fulfills the condition above, and so we will call it the $N$-truncated discrete generalized exchange-driven system 
($N$-DGED for short):
\begin{equation}
\left\{
\begin{aligned}
\dot{c}_0  &= \sum_{j=1}^{2} Q_{j, 0}^N(c)    \\
\dot{c}_i  &= \sum_{j=1}^{4} Q_{j, i}^N(c), \quad i\in \{1, \ldots, N-1\}   \\
\dot{c}_N  &= \sum_{j=3}^{4} Q_{j, N}^N(c),    
\end{aligned}
\right. \label{NDGED}
\end{equation}
where
\begin{align}
Q_{1, i}^N(c) & := \sum_{k=1}^{N-i} \sum_{j=0}^{N-k} a(i+k, j; k)c_{i+k}c_{j},		\label{Q1Ni} \\
Q_{2, i}^N(c) & := - \sum_{k=1}^{N-i} \sum_{j=k}^N a(j, i; k)c_{j}c_{i},		\label{Q2Ni} \\
Q_{3, i}^N(c) & := \sum_{k=1}^i \sum_{j=k}^N a(j, i-k; k)c_{j}c_{i-k},		\label{Q3Ni} \\
Q_{4, i}^N(c) & := -  \sum_{k=1}^i \sum_{j=0}^{N-k} a(i, j; k)c_{j}c_{i},		\label{Q4Ni} 
\end{align}
and
\begin{equation} \label{Tinit}
c_i^N(0) = c_{0 i}\geq 0,\qquad i\in\{0,\ldots,N\}.
\end{equation}

\medskip

In the non-isolated case the $N$-DGED system is the same as in \eqref{NDGED} 
but for the $c_0$-equation which is substituted by $\dot{c}_0=0,$
or, equivalently, we define, for all $j$ and $N$, 
\begin{equation}
    Q_{j,0}^N =  0. \label{Qj0N}
\end{equation}

\medskip

Being an ordinary differential equation in $\Rb^{N+1}$ with a polynomial vector field, the existence and uniqueness of solutions
to the initial value problem are ensured by the standard Picard-Lindel\"of theorem
(see, e.g.,  \cite[Theorem I-1-4]{HS}). It is also easy to conclude by
standard arguments (see, e.g., the proof of Theorem III-4-5 in \cite{HS}) that $\Rb^{(N+1)+}$ is invariant for the local flow 
associated with \eqref{NDGED}, which means that nonnegative initial data have unique nonnegative local solutions. 
The following result is important for the remaining
analysis

\begin{prop}\label{Prop1}
Let $c^N=(c^N_i)_{0 \leq i \leq N }$ be any solution of \eqref{NDGED}--\eqref{Q4Ni} in the isolated case. 
Then, for every $(g_i)$ we have
\begin{equation}
\frac{d}{dt}\sum_{i=0}^Ng_ic_i^N = \sum_{k=1}^{N-1}\sum_{i=k}^N\sum_{j=0}^{N-k}(g_{j+k}+g_{i-k}-g_j-g_i)a(i, j; k)c_i^Nc_j^N \label{gNmom}
\end{equation}
\end{prop}
\begin{Proof}
Let us first take the summation of the quantity $g_i \dot{c}_i^N $ from $i=0$ to $i=N$. Then, from \eqref{NDGED}, it can be inferred that
\begin{equation}\label{PropEqua1}
\begin{split}
\frac{d}{dt}\sum_{i=0}^N &g_i c_i^N 
 =  g_0 \sum_{k=1}^N \sum_{j=0}^{N-k}  a(k, j; k) c_k^N c_j^N  - g_0 \sum_{k=1}^N \sum_{j=k}^{N}  a(j, 0; k) c_j^N c_0^N \\
 & + \sum_{i=1}^{N-1} \sum_{k=1}^{N-i} \sum_{j=0}^{N-k} g_i a(i+k, j; k) c_{i+k}^N c_j^N -  \sum_{i=1}^{N-1} \sum_{k=1}^{N-i} \sum_{j=k}^{N} g_i a(j, i; k) c_{j}^N c_i^N \\
 & + \sum_{i=1}^{N-1} \sum_{k=1}^{i} \sum_{j=k}^{N} g_i a(j, i-k; k) c_{j}^N c_{i-k}^N -  \sum_{i=1}^{N-1} \sum_{k=1}^{i} \sum_{j=0}^{N-k} g_i a(i, j; k) c_{i}^N c_j^N \\
 & +  g_N \sum_{k=1}^{N} \sum_{j=k}^{N}  a(j, N-k; k) c_{j}^N c_{N-k}^N - g_N \sum_{k=1}^{N} \sum_{j=0}^{N-k}   a(N, j; k) c_{j}^N c_N^N.
\end{split}
\end{equation}
Next, we simplify some terms on the right-hand side of \eqref{PropEqua1}. To do this, we will start with the third sum on the right-hand side of \eqref{PropEqua1}. Changing the order of summation of the sums in $i$ and $k$ and renaming $i+k\mapsto i$, it can be rewritten 
\begin{align}\label{PropEqua2}
  \sum_{i=1}^{N-1} \sum_{k=1}^{N-i} \sum_{j=0}^{N-k} g_i a(i+k, j; k) c_{i+k}^N c_j^N = & \sum_{k=1}^{N-1} \sum_{i=1}^{N-k} \sum_{j=0}^{N-k} g_i a(i+k, j; k) c_{i+k}^N c_j^N  \nonumber\\
= & \sum_{k=1}^{N-1} \sum_{i=k+1}^{N} \sum_{j=0}^{N-k} g_{i-k} a(i, j; k) c_i^N c_j^N  \nonumber\\
= & \sum_{k=1}^{N-1} \sum_{i=k}^{N} \sum_{j=0}^{N-k} g_{i-k} a(i, j; k) c_i^N c_j^N  \nonumber\\
 & - g_{0}  \sum_{k=1}^{N-1}  \sum_{j=0}^{N-k}  a(k, j; k) c_k^N c_j^N.
\end{align}
In the same vein, the fourth sum in \eqref{PropEqua1} 
can be rearranged by changing the order of the sums twice, first to exchange the $i$ and $k$ sums and then 
the $i$ and $j$ sums, and finally changing notation $i\leftrightarrow j,$ giving
\begin{multline}\label{PropEqua3}
    \sum_{i=1}^{N-1} \sum_{k=1}^{N-i} \sum_{j=k}^{N} g_i a(j, i; k) c_{j}^N c_i^N\\
    =  \sum_{k=1}^{N-1}  \sum_{i=k}^{N} \sum_{j=0}^{N-k} g_j a(i, j; k) c_{j}^N c_i^N - g_0 \sum_{k=1}^{N-1}  \sum_{i=k}^{N}   a(i, 0; k) c_{0}^N c_i^N. 
\end{multline}
Similarly, the fifth sum can be rearranged by exchanging the $i$ and $k$ sums,
then introducing $i'=i-k,$ changing again the order of the $i'$ and $j$ sums, 
and finally change notation $j\mapsto i$ and $i'\mapsto j$. In the end we get
\begin{multline}\label{PropEqua4}
 \sum_{i=1}^{N-1} \sum_{k=1}^{i} \sum_{j=k}^{N} g_i a(j, i-k; k) c_{j}^N c_{i-k}^N  \\
=  \sum_{k=1}^{N-1}  \sum_{i=k}^{N}  \sum_{j=0}^{N-k}  g_{j+k} a(i, j; k) c_{i}^N c_{j}^N 
 -  g_{N} \sum_{k=1}^{N-1}  \sum_{i=k}^{N}    a(i, N-k; k) c_{i}^N c_{N-k}^N.
\end{multline}
Furthermore, let us write the sixth summation in simplified form by changing the order of the $i$ and $k$ sums:
\begin{multline}\label{PropEqua5}
\sum_{i=1}^{N-1} \sum_{k=1}^{i} \sum_{j=0}^{N-k} g_i a(i, j; k) c_{i}^N c_j^N  \\
=  \sum_{k=1}^{N-1} \sum_{i=k}^{N} \sum_{j=0}^{N-k} g_i a(i, j; k) c_{i}^N c_j^N - \sum_{k=1}^{N-1}  \sum_{j=0}^{N-k} g_N a(N, j; k) c_{N}^N c_j^N.
\end{multline}
Let us now incorporate the results obtained from equations \eqref{PropEqua2}--\eqref{PropEqua5} into equation \eqref{PropEqua1}. The revised expression after these substitutions is
\begin{equation}\label{PropEqua6}
\begin{split}
 \frac{d}{dt} \sum_{i=0}^N &g_i c_i^N 
=   \sum_{k=1}^{N-1} \sum_{i=k}^{N} \sum_{j=0}^{N-k} [ g_{i-k} + g_{j+k} -g_i - g_j ]  a(i, j; k) c_i^N c_j^N  \\ 
& +  g_0 \sum_{k=1}^N \sum_{j=0}^{N-k}  a(k, j; k) c_k^N c_j^N   - g_0 \sum_{k=1}^N \sum_{j=k}^{N}  a(j, 0; k) c_j^N c_0^N  \\
  &  -  g_{0} \sum_{k=1}^{N-1}  \sum_{j=0}^{N-k}  a(k, j; k) c_k^N c_j^N   +  g_0 \sum_{k=1}^{N-1}  \sum_{i=k}^{N}   a(i, 0; k) c_{0}^N c_i^N  \\ 
 &  -  \sum_{k=1}^{N-1}  \sum_{i=k}^{N}   g_{N} a(i, N-k; k) c_{i}^N c_{N-k}^N  + \sum_{k=1}^{N-1}  \sum_{j=0}^{N-k} g_N a(N, j; k) c_{N}^N c_j^N  \\
&  + \sum_{k=1}^{N} \sum_{j=k}^{N} g_N a(j, N-k; k) c_{j}^N c_{N-k}^N  
 -  \sum_{k=1}^{N} \sum_{j=0}^{N-k}  g_N a(N, j; k) c_{j}^N c_N^N.
\end{split}
\end{equation}
Furthermore, we again write \eqref{PropEqua6} into the following concise form:
\begin{align*}
 \frac{d}{dt} \sum_{i=0}^N g_i c_i^N 
= &  \sum_{k=1}^{N-1} \sum_{i=k}^{N} \sum_{j=0}^{N-k} [ g_{i-k} + g_{j+k} -g_i - g_j ]  a(i, j; k) c_i^N c_j^N  \nonumber\\  & +  g_0   a(N, 0; N) c_N^N c_0^N   - g_0  a(N, 0; N) c_N^N c_0^N     \nonumber\\ 
 & +   g_N a(N, 0; N) c_{N}^N c_{0}^N  - g_N a(N, 0; N) c_{0}^N c_N^N \nonumber\\
= &  \sum_{k=1}^{N-1} \sum_{i=k}^{N} \sum_{j=0}^{N-k} [ g_{i-k} + g_{j+k} -g_i - g_j ]  a(i, j; k) c_i^N c_j^N.
\end{align*}
This completes the proof of Proposition \ref{Prop1}.
\end{Proof}

\medskip

It is easy to observe that the right-hand side of \eqref{gNmom} is identically zero when $g_i=i$, and $g_i=1$, thus proving the following
two conservation laws:
\begin{corol}\label{Corollary1}
All non negative solutions to \eqref{NDGED} in the isolated case conserve the total number of initial clusters $\|c(0)\|_{\ell_1}$ and the initial mass $\|(ic_i(0))\|_{\ell_1},$
when we consider a solution $c^N$ of \eqref{NDGED} an $\ell_1$-valued function by defining $c_j^N\equiv 0$ for all $j\geq N+1.$
\end{corol}

\begin{remark} A consequence of the previous corollary is that any non negative solution of the truncated system \eqref{NDGED} is globally defined in $\Rb^+_0$.\end{remark}

\medskip

In the non-isolated case there is a similar result whose proof proceed in the same way and will be omitted. 
\begin{prop}\label{Prop1bis}
Let $c^N=(c^N_i)_{0 \leq i \leq N }$ be any solution of \eqref{NDGED}--\eqref{Tinit} in the non-isolated case. 
Then, for every $(g_i)$ we have
\begin{equation}\label{gNmomnonisolated}
\begin{split}
\frac{d}{dt}\sum_{i=0}^Ng_ic_i^N  & = \sum_{k=1}^{N-1}\sum_{i=k}^N\sum_{j=0}^{N-k}(g_{j+k}+g_{i-k}-g_j-g_i)a(i, j; k)c_i^Nc_j^N  \\
& \quad-g_0\Bigl(\sum_{k=1}^N\sum_{j=0}^{N-k}a(k,j;k)c_k^Nc_j^N - \sum_{k=1}^N\sum_{i=k}^Na(i,0;k)c_i^Nc_{0}^N\Bigr).
\end{split}
\end{equation}
\end{prop}

\medskip

\begin{remark} Taking $g_i=i$ in \eqref{gNmomnonisolated}
it is easily seen that solutions of the DGED system in the 
 non-isolated case conserve the initial mass and, as
 a consequence, any non negative solution of the truncated system \eqref{NDGED} is globally defined in $\Rb^+_0$
 also in the non-isolated case. However,
 making $g_i=1$ in \eqref{gNmomnonisolated} we immediately 
 conclude that the total number of initial clusters is no longer conserved in this case.\end{remark}

\medskip

The following two lemmas about properties of the 
solutions to the truncated system, in which we assume the conditions \eqref{qikbound}-\eqref{bound},
will be essential for the proof of the existence theorem in the next section. They use, in a
crucial way, the inequalities proved in lemmas~\ref{sigmalemma} and~\ref{sigmalemma2}. 

\begin{lemma}\label{HMomentLem1}
Let $c^N = (c_i^N)_{ 0 \leq i \leq N }$ be a solution to \eqref{NDGED}--\eqref{Tinit} and $\sigma  \in \mathcal{E} \cup \mathcal{E}_1$. Assume $ \sum_{i=0}^N \sigma(i) c_{0 i}^N$ is finite. Then, for each $T\geq 0,$ we have for all $t\in[0,T],$
\begin{equation}\label{gamma1lemma}
 \sum_{i=0}^N \sigma(i) \  c_i^N (t) \leq \gamma_T \sum_{i=0}^N \sigma(i) c_{0 i}^N,
\end{equation}
and, for integers $M,N$ such that $M_0\leq M\leq N$, with $M_0$ as in Lemma \ref{sigmalemma2},
\begin{equation}\label{gamma1lemma1}
0\leq \int_0^T \sum_{(p,k)\in\mathcal{J}^N_0} \frac{\sigma(p-1)}{p-1}a(p,0;k)c^N_0(t)c^N_p(t)dt
\leq \gamma_T \sum_{i=0}^N \sigma(i) c_{0 i}^N,
\end{equation}
 where, 
 $$
 \mathcal{J}^N_0:=\{(p,k)\in(\Nb_+)^2:M_0\leq p\leq N,\; 1\leq k\leq p-1\},
 $$
 and $\gamma_T>0$ is a constant only depending on $C, m_\sigma,\mathcal{P}_1(0),\mathcal{Q}, T. $ 
\end{lemma}
\begin{proof}
For $1\leq i \leq N$, it can be deduced from Proposition \ref{Prop1}, by setting $g_i := \sigma(i), $ that
 \begin{equation}\label{DLarge1}
\frac{d}{dt}\sum_{i=0}^N \sigma(i) c_i^N(t) = \sum_{k=1}^{N-1}\sum_{i=k}^N\sum_{j=0}^{N-k} \tilde{ \sigma } (i,j,k) \ a(i, j; k) \  c_i^N(t) c_j^N(t),
\end{equation}
where
\begin{equation}
 \tilde{ \sigma }(i,j,k)  :=   \sigma(j+k) + \sigma(i-k) - \sigma(j) - \sigma(i).\label{sigmatest}
\end{equation}
Since there are no upper bounds on the rate coefficients of the type $a(p,0;k),$ the terms involving this type of coefficients have to be tackled differently from the others.

For $j\in\Nb_+$, by  \eqref{bound}, \eqref{sigmatest} and Lemma \ref{sigmalemma}, we obtain,
\begin{align}\label{DLarge2}
\tilde{ \sigma } (i, j, k)    a(i, j; k) \leq  & C(i-k+1)(j+k)q_{i,k}\tilde{\sigma}(i,j,k) \nonumber\\
\leq  & Cm_\sigma( j \sigma(k) + k \sigma(j) )(i-k+1)q_{i,k},
\end{align}
so, by applying \eqref{DLarge2} to  the terms not involving the rate coefficients  of the type $a(p,0;k),$ we have, for $N\geq 2,$
 \begin{align}\label{S0est}
\sum_{k=1}^{N-1}\sum_{i=k}^N\sum_{j=1}^{N-k} &\tilde{ \sigma } (i,j,k) \ a(i, j; k) \  c_i^N(t) c_j^N(t)\\ \leq & Cm_\sigma \left(\sum_{k=1}^{N-1}\sum_{i=k}^N\sum_{j=1}^{N-k} j\sigma(k)(i-k+1) q_{i,k}c_i^N(t) c_j^N(t)\right.\\
&\qquad\qquad\qquad\left. + 
\sum_{k=1}^{N-1}\sum_{i=k}^N\sum_{j=1}^{N-k} k\sigma(j) (i-k+1)q_{i,k} c_i^N(t) c_j^N(t)\right) \nonumber\\
{\color{blue}=:} & Cm_\sigma(S_1+S_2).
\end{align}

We now estimate separately $S_1$ and $S_2$. By using the convexity of $\sigma$, which entails the fact that $k\mapsto \sigma(k)/k$ is increasing, hypothesis \eqref{qikbound} and
Corollary \ref{Corollary1}, we obtain, 
\begin{align}\label{S1est}
S_1&=\sum_{k=1}^{N-1}\sum_{i=k}^N\sum_{j=1}^{N-k}\frac{\sigma(k)}{k}\frac{i}{\sigma(i)}\left(\frac{k(i-k+1)}{i}q_{i,k}\right) ( j  c_j^N(t))(\sigma(i) c_i^N(t)) \nonumber\\
&\leq  \mathcal{P}_1(0)\mathcal{Q}  \sum_{i=0}^{N}\sigma(i) c_i^N(t).
\end{align}
Also,
\begin{align}\label{S2est}
S_2&=\sum_{k=1}^{N-1}\sum_{i=k}^N\sum_{j=1}^{N-k}\left(\frac{k(i-k+1)}{i}q_{i,k}\right) (ic_i^N(t))(\sigma(j) c_j^N(t))\nonumber\\
&\leq   \mathcal{P}_1(0)\mathcal{Q}  \sum_{j=0}^{N}\sigma(j) c_j^N(t).
\end{align}
By applying \eqref{S0est}--\eqref{S2est} to \eqref{DLarge1} we obtain
\begin{equation}\label{S3est}
\frac{d}{dt}\sum_{i=0}^N \sigma(i) c_i^N(t)\leq K\sum_{i=0}^N \sigma(i)c_i^N(t)
+\sum_{k=1}^{N-1}\sum_{i=k}^{N}\tilde{\sigma}(i,0,k)a(i,0;k)c_0^N(t)c_i^N(t).
\end{equation}
where $K:=2C\mathcal{Q}m_\sigma\mathcal{P}_1(0).$ But, according to the proof of Lemma \ref{sigmalemma}, for $1\leq k\leq i,$
$$
\tilde{\sigma}(i,0,k)=-(\sigma(i)-\sigma(i-k)-\sigma(k))\leq 0,
$$
from \eqref{S3est} we obtain
\begin{equation}\label{S4est}
\frac{d}{dt}\sum_{i=0}^N \sigma(i) c_i^N(t)\leq K\sum_{i=0}^N \sigma(i)c_i^N(t),
\end{equation}
thus resulting \eqref{gamma1lemma} by Gronwall inequality. 

In order to prove \eqref{gamma1lemma1}, we integrate both members of \eqref{S3est} and take {\color{blue} into} account \eqref{gamma1lemma} and the sign of each term, thus getting
\begin{align}\label{L1s5est}
\int_0^T\Big|\sum_{k=1}^{N-1}\sum_{i=k}^{N}\tilde{\sigma}(i,0,k)a(i,0;k)c_0^N(t)c_i^N(t)\Big|dt
\leq&\sum_{i=0}^N\sigma(i)c^N_i(0)+K\int_0^T\sum_{i=0}^N \sigma(i) c_i^N(t)dt\nonumber\\
\leq&\bar{\gamma}_T \sum_{i=0}^N\sigma(i)c^N_i(0),
\end{align}
where $\bar{\gamma}_T:=1+KT\exp\big(KT\big).$ Remarking that $\tilde{\sigma}(i,0,i)=0,$ we can write, for each $i=1,\ldots,j,$ and $N\geq 2,$
\begin{multline}\label{L1s6est}
\int_0^T\Big|\sum_{k=1}^{N-1}\sum_{i=k}^{N}\tilde{\sigma}(i,0,k)a(i,0;k)c_0^N(t)c_i^N(t)\Big|dt\\= 
\int_0^T\Big|\sum_{i=1}^{N}\sum_{k=1}^{i-1}\tilde{\sigma}(i,0,k)a(i,0;k)c_0^N(t)c_i^N(t)\Big|dt.
\end{multline}
Since $-\tilde{\sigma}(i,0,k)=|\tilde{\sigma}(i,0,k)|=\sigma(i)-\sigma(i-k)-\sigma(k),$ the conclusion is obtained from \eqref{L1s5est} and \eqref{L1s6est} by application of Lemma \ref{sigmalemma2} and redefining the constant $\gamma_T.$ 


\end{proof}


\begin{lemma}\label{tderboundlemma}
Let $c^N = (c_i^N)_{ 0 \leq i \leq N }$ be a solution to \eqref{NDGED} and $T\in(0,+\infty).$ There exists a positive constant $ \Gamma_{2,i}(T)$ depending only 
on $C, \mathcal{Q}, \|c_0\|, i$ and $T$ such that, for each $i\in\Nb,$
\begin{equation*}
\bigg| \frac{dc_i^N }{dt} \bigg|_{L^{1}(0, T) } \leq\Gamma_{2,i}(T).
\end{equation*}
\end{lemma}

\begin{proof}
We present the proof for the isolated case and will
comment on the non-isolated case at the end.

\medskip

By Corollary \ref{Corollary1}, the following expressions can be deduced:
\begin{align}\label{TMassEst}
    \sum_{i=0}^{N} i c_i^N(t) =     \sum_{i=0}^{N} i c_{0 i}^N \leq   \sum_{i=0}^{\infty} i c_{0 i}^N =\| (i c_{i0} )\|_{\ell_1}\leq \|c_0\|,
\end{align}
and then
\begin{align}\label{TParticleEst}
    \sum_{i=0}^{N}  c_i^N(t) \leq \sum_{i=0}^{N}  ic_i^N(t) \leq\|c_0\|.
\end{align}
Also recall that the numbers $q_{i,k}$ satisfy estimate \eqref{qikest}.
\medskip

Now, consider system \eqref{NDGED}. Let $\hat{Q}^N_{j,i},$ $j=1,2,3,4,$ be the sums \eqref{Q1Ni}-\eqref{Q4Ni} taken away the terms of the type $a(p,0;k)c^N_{p}c^N_0$. Hence, considering that $\hat{Q}^N_{2,0}=\hat{Q}^N_{3,1}=0,$ we can write
\begin{align}\label{hatQN1I}
Q^N_{1,i}(c^N)&=\hat{Q}^N_{1,i}(c^N)+\sum_{k=1}^{N-i}a(i+k,0;k)c^N_{i+k}c^N_0,\;\;\quad\qquad  i\in\{0,\dots,N-1\},\\
\label{hatQN2I}Q^N_{2,i}(c^N)&=(1-\delta_{i,0})\hat{Q}^N_{2,i}(c^N) - \delta_{i,0}{\displaystyle \sum_{k=1}^{N}\sum_{j=k}^Na(j,0;k)c^N_jc^N_0,} \;\;\; i\in\{0,\dots,N-1\},\\
\label{hatQN3I}Q^N_{3,i}(c^N)&={\displaystyle (1-\delta_{i,1})  \hat{Q}^N_{3,i}(c^N)+\displaystyle\sum_{j=i}^Na(j,0;i)c^N_jc^N_0,}\quad\;\; \quad i\in\{1,\dots,N\},\\
\label{hatQN4I}Q^N_{4,i}(c^N)&=\hat{Q}^N_{4,i}(c^N)-\sum_{k=1}^ia(i,0;k)c^N_0c^N_i\,,\qquad\qquad\;\;\; \quad i\in\{1,\dots,N\},
\end{align}
where $\delta_{i,\ell}$ is the Kronecker symbol.
Therefore, after some rearrangements, we can rewrite the truncated system \eqref{NDGED} in the following form:
\begin{align}
\dot{c}^N_0  &= \hat{Q}^N_{1,0}(c^N)-\sum_{k=1}^{N}\sum_{j=k+1}^Na(j,0;k)c^N_jc^N_0, \label{NDGED10} \\
\dot{c}^N_i  &= \sum_{j=1}^{4} \hat{Q}^N_{j, i}(c)+2\sum_{j=i+1}^N a(j,0;i)c^N_jc^N_0\nonumber\\
&\qquad\qquad\qquad\qquad-\sum_{k=1}^{i-1}a(i,0;k)c^N_ic^N_0,\qquad\text{for }i\in\{1,\dots,N-1\}, \label{NDGED1i}  \\
\dot{c}^N_N  &= \sum_{j=3}^{4} \hat{Q}_{j, N}^N(c)-\sum_{k=1}^{N-1}a(N,0;k)c^N_Nc^N_0.   \label{NDGED1N} 
\end{align}

Since no upper bounds are imposed on the fragmentation-type coefficients $a(p,0;k),$ the terms involving these must be estimated in a separate way, as in \cite{PL:2002}.
\medskip

We first proceed to the estimation of the $\hat{Q}^N_{j,i}$ terms.  We remark that we are using in this work the convention that if the lower index of a sum is greater than the upper index then the sum is zero.  Taking into account \eqref{Q1Ni}, we have, for $i=0,\dots,N-1,$
\begin{align}\label{Q1NiA}
\hat{Q}^N_{1,i}(c^N)&\leq C\sum_{k=1}^{N-i}\sum_{j=1}^{N-k}(i+1)(j+k)q_{i+k,k}c^N_{i+k}c^N_j\nonumber\\
&\leq C\mathcal{Q}\sum_{k=1}^{N-i}\sum_{j=1}^{N-k}\big((i+k)c^N_{i+k}\big)(jc^N_j)
+C\sum_{k=1}^{N-i}\sum_{j=1}^{N-k}(i+1)kq_{i+k,k}c^N_{i+k}c^N_j\nonumber\\
&\leq C\mathcal{Q}\|(ic_{i0})\|_{\ell_1}^2+C\mathcal{Q}(i+1)\|(ic_{i0})\|_{\ell_1}\|c_0\|_{\ell_1}\nonumber\\
&\leq 2C\mathcal{Q}(i+1)\|c_0\|^2,
\end{align}
where we have used \eqref{qikest} which, in this case, implies that $kq_{i+k,k}\leq \mathcal{Q}(i+k).$

By \eqref{Q2Ni} we have, for $i=1,\dots,N-1,$
\begin{align}\label{Q2NiA}
|\hat{Q}^N_{2,i}(c^N)|&\leq C\sum_{k=1}^{N-i}\sum_{j=k}^{N}(j-k+1)(i+k)q_{j,k}c^N_{j}c^N_i\nonumber\\
&\leq C\sum_{j=1}^{N}\sum_{k=1}^{j}(j-k+1)q_{j,k}c^N_{j}(ic^N_i)
+C\sum_{j=1}^{N}\sum_{k=1}^{j}(j-k+1)kq_{j,k}c^N_{j}c^N_i\nonumber\\
&\leq C\mathcal{Q}\|(ic_{i0})\|_{\ell_1}^2+C\mathcal{Q}\|(ic_{i0})\|_{\ell_1}\|c_0\|_{\ell_1}
\leq 2C\mathcal{Q}\|c_0\|^2.
\end{align}
where again we have used \eqref{qikbound} that implies this time that $(j-k+1)kq_{j,k}\leq\mathcal{Q}j.$

\medskip

By \eqref{Q3Ni}, for $i=2,\dots,N,$
\begin{align}\label{Q3NiA}
\hat{Q}^N_{3,i}(c^N)&\leq C\sum_{k=1}^{i-1}\sum_{j=k}^{N}(j-k+1)iq_{j,k}c^N_{j}c^N_{i-k}\nonumber\\
&\leq C\mathcal{Q}i\|(ic_{i0})\|_{\ell_1}\|(c_{i0})\|_{\ell_1}\leq C\mathcal{Q}i\|c_0\|^2.
\end{align}
Finally, from \eqref{Q3Ni}, for each $i=1,\dots,N$  we obtain
\begin{align}\label{Q4NiA}
|\hat{Q}^N_{4,i}(c^N)|
&\leq C\sum_{k=1}^i \sum_{j=1}^{N-k} (i-k+1)(j+k)q_{i,k}c^N_{j}c^N_{i}\nonumber\\
&\leq C\sum_{k=1}^i \sum_{j=1}^{N-k} (i-k+1)q_{i,k}(jc^N_{j})c^N_{i}+
 C\sum_{k=1}^i \sum_{j=1}^{N-k} (i-k+1)q_{i,k}kc^N_{j}c^N_{i}\nonumber\\
&\leq 2C\mathcal{Q}\|c_0\|^2, 
\end{align}
where, once more, we have used \eqref{qikbound} in a way similar to the cases above.

\medskip

To proceed to the estimates on the other terms, in the case $i=1,\dots,N-1,$ 
we first integrate \eqref{NDGED1i} to get
\begin{multline*}
c^N_i(T)-c^N_i(0)=\int_0^T\sum_{j=1}^{4}\hat{Q}^N_{j,i}(c^N(t))dt\\
+2\int_0^T\sum_{j=i+1}^N a(j,0,;i)c^N_j(t)c^N_0(t)dt-\int_0^T \sum_{k=1}^{i}a(i,0;k)c^N_i(t)c^N_0(t)dt,
\end{multline*}
where we take $\hat{Q}^N_{3,1}\equiv 0.$
Considering the signs of the $\hat{Q}^N_{j,i},$  we can infer that
\begin{align}\label{estint1}
0&\leq 2\int_0^T\sum_{j=i+1}^N a(j,0,;i)c^N_j(t)c^N_0(t)dt\nonumber\\
&\leq c^N_i(T)+\int_0^T \sum_{j\in\{2,4\}}|\hat{Q}^N_{j,i}(c^N(t))|dt+\int_0^T \sum_{k=1}^{i}a(i,0;k)c^N_i(t)c^N_0(t)dt.
\end{align}
On the other hand, calling $\alpha(i):=\sum_{k=1}^{i}a(i,0;k)$ we obtain, for the last term in \eqref{estint1},
\begin{equation}\label{estint2}
\int_0^T \sum_{k=1}^{i}a(i,0;k)c^N_i(t)c^N_0(t)dt\leq \alpha(i)\|c_0\|^2 T.
\end{equation}
Therefore, taking into account \eqref{TParticleEst}, \eqref{Q2NiA}, \eqref{Q4NiA}, \eqref{estint1} and \eqref{estint2},  we obtain
\begin{equation}\label{L1estA}
\int_0^T\sum_{j=i+1}^N a(j,0;i)c^N_j(t)c^N_0(t)dt\leq \gamma.
\end{equation}
for $\gamma:= \frac{1}{2}\|(c_{i0})\|_{\ell_1}+
2C\mathcal{Q}\|(ic_{i0})\|_{\ell_1}^2T+\frac{\alpha(i)}{2}\|(c_{i0})\|_{\ell_1}^2T.$
Then, using estimates \eqref{estint2} and \eqref{L1estA},  together with \eqref{Q1NiA}-\eqref{Q4NiA} in the equation \eqref{NDGED1i}, we conclude that, for
$i=1,\ldots,N-1,$
$$
|\dot{c}^N_i|_{L^1(0,T)}\leq
\bigl(9C\mathcal{Q}(i+1)+\frac{3}{2}\alpha(i)\bigr)\|c_0\|^2T+\frac{1}{2}\|c_0\|,
$$
which proves our claim for $i=1,\ldots,N-1.$ We omit the proof for the cases $i=0$ and $i=N,$ since they are easier and follow exactly the same lines. For $i\geq N+1$ it is trivial.

\medskip

For the non-isolated case the only difference is that 
all $Q_{j,0}^N$ are identically zero, instead of \eqref{hatQN1I} and
\eqref{hatQN2I} with $i=0$,
and so the estimates above still hold.
\end{proof}

%
%
\section{Existence  for the general exchange-driven system}\label{sec3}

\subsection{The isolated case \eqref{DGED}}
In this section 
prove the existence of solutions for the initial value problem \eqref{DGED}, \eqref{init}
in the isolated system case.



\begin{Proofof} {\color{blue}of theorem \ref{ExistThm1}:} By using a refined version of de la Vallée-Poussin theorem as in  \cite[Section 4]{PL:2002}, we can conclude that our hypothesis for $c_0$ implies that there is $\sigma_0\in\mathcal{E}$ such that 
\begin{equation}\label{sigmac0}
\mathcal{S}_0:=\sum_{i=0}^\infty \sigma_0(i)c_{i0}<\infty.
\end{equation}
Then, since for any $N\in\Nb_+,$ 
$$
\sum_{i=0}^N \sigma_0(i)c^N_{i0}=\sum_{i=0}^N \sigma_0(i)c_{i0}\leq \mathcal{S}_0,
$$ 
we conclude that, for each $T\geq 0,$ we have, for all $t\in[0,T],$
\begin{equation}\label{sigmact}
\sum_{i=0}^N \sigma_0(i)c^N_{i}(t)\leq \gamma_T \mathcal{S}_0,
\end{equation}
by using Lemma \ref{HMomentLem1}.

It follows from Lemma~\ref{tderboundlemma} that, for each $i\in\Nb,$ $ (c^N_i )_{N\geq i}$  is a bounded sequence in $L^\infty(0,T)\cap W^{1,1}(0,T)$. Therefore, we can apply the Helly selection theorem, together with a diagonal argument, to guarantee the existence of a subsequence  $(c^{N}_i)_{N\geq i}$ (not relabelled) and a  sequence $c=(c_i)_{i\in\Nb}$ of nonnegative functions of locally bounded variation such that, for each $i\in\Nb$ and  $t \in [0, \infty )$,
\begin{equation}\label{Converge1}
\lim_{N\to +\infty}c^N_i(t) = c_i(t).
\end{equation}
Also, it follows from the nonnegativity of the $c^N_i$, \eqref{TMassEst}, \eqref{TParticleEst}, and \eqref{Converge1}, that, for each $t\in[0,+\infty),$ we have $c(t)\in X_{0,1}^+$ and
\begin{equation}\label{normest}
\|c(t)\|\leq \|c_0\|.
\end{equation}
Let $T>0.$ Fixing any integer $M\geq 1$, adding in \eqref{sigmact} up to $M$ and then taking limit as $N\to +\infty$ in $c^N_i(t)$ in this sum, 
we obtain the same inequality with $c^N_i$ replaced by $c_i$. Then, taking $M\to +\infty$, we infer that
\begin{equation}\label{sigmactlim}
\sum_{i=0}^\infty \sigma_0(i)c_{i}(t)\leq \gamma_T\mathcal{S}_0,\qquad t\in [0,T].
\end{equation}

\medskip

Now we prove that $c$ is indeed a solution to \eqref{DGED}.  According to Definition~\ref{Def Sol} it is necessary that, for $j\in\{1,2,3,4\}$ and $i\in\Nb_+,$ $Q_{j,i}(c(\cdot))\in L^1(0,T).$ Similarly to the proof of Lemma \ref{tderboundlemma}, we define the quantities $\hat{Q}_{j,i}(c),$ $i\in\Nb,$ $j\in\{1,2,3,4\},$ as the corresponding $Q_{j,i}(c)$ taken away the terms of the type $a(p,0;k)c_pc_0.$ 
We recall that this procedure is forced upon us by the fact that
the growth bound \eqref{bound} on the rate coefficients 
$a(i,j;k)$ does not necessarily hold when $j=0$.

\medskip

We claim that, for each $i\in\Nb_+$ and $j\in\{1,2,3,4\},$ 
\begin{equation}\label{claim1}
\hat{Q}_{j,i}(c)\in L^1(0,T).
\end{equation} 
Define first
$$
\Phi(\beta):=\sum_{(i,j,k)\in \mathcal{I}}a(i,j;k)\beta_i\beta_j
$$
where $\mathcal{I}:=\{(i,j,k)\in(\Nb_+)^3:k\leq i\},$ with domain $D_\Phi$ made by the sequences $\beta\in X_{0,1}^+$ for which the above triple series is convergent. We prove that, in fact, $D_\Phi=X_{0,1}^+.$ Let $\eta\in\Nb_+.$ Then, for each $\beta\in X_{0,1}^+,$ by \eqref{qikbound} and \eqref{bound}, we have,
\begin{align*}
\sum_{i=1}^\eta\sum_{j=1}^\eta\sum_{k=1}^ia(i,j;k)\beta_i\beta_j\leq\; & 
C\sum_{i=1}^\eta\sum_{j=1}^\eta\sum_{k=1}^i (i-k+1)(j+k)q_{i,k}\beta_i\beta_j\\
=\; & C\sum_{i=1}^\eta\sum_{j=1}^\eta\sum_{k=1}^i (i-k+1)q_{i,k}\beta_i(j\beta_j)\\&+
C\sum_{i=1}^\eta\sum_{j=1}^\eta\sum_{k=1}^i (i-k+1)kq_{i,k}\beta_i\beta_j\\
\leq\; & C\mathcal{Q}\Big(\sum_{i=1}^\eta i\beta_i\Big)\Big(\sum_{j=1}^\eta j\beta_j\Big)+
C\mathcal{Q}\Big(\sum_{i=1}^\eta i\beta_i\Big)\Big(\sum_{j=1}^\eta \beta_j\Big)\\
\leq \;& 2C\mathcal{Q}\|\beta\|^2.
\end{align*}
By making $\eta\to \infty,$ we conclude that the above triple series is convergent and therefore, $\beta \in D_\Phi.$ Now observe that, for each $j\in\{1,2,3,4\}$ and $i\in\Nb_+,$ the indices in the double sum in $\hat{Q}_{j,i},$ runs over a subset of $\mathcal{I},$ so that, for each $\beta\in X_{0,1}^+,$ $\hat{Q}_{j,i}(\beta)$ is well defined and
\begin{equation*}
|\hat{Q}_{j,i}(\beta)|\leq \Phi(\beta)\leq 2C\mathcal{Q}\|\beta\|^2.
\end{equation*}
Hence, by \eqref{normest}, 
\begin{equation}\label{hatQPhibound}
\int_0^T |\hat{Q}_{j,i}(c(t))|dt\leq  2C\mathcal{Q}\|c_0\|^2 T,
\end{equation}
and our claim \eqref{claim1} is proved.


Our next claim is that, for $j\in\{1,2,3,4\},$ $i\in\Nb_+,$ and also for $(j,i)=(1,0),$
\begin{equation}\label{claim2}
\lim_{N\to +\infty}|\hat{Q}_{j,i}^N(c^N)-\hat{Q}_{j,i}(c)|_{L^1(0,T)}=0.
\end{equation}
We first rearrange the sums that define $\hat{Q}_{j,i}^N$ and rewrite them, for each sequence $\beta$ as 
$$
\hat{Q}_{j,i}^N(\beta)=\sum_{(p,q)\in\mathcal{J}^N_{j,i}}\;\sum_{k\in\mathcal{K}^N_{j,i}(p,q)}a(p,q;k)\beta_p\beta_q,
$$
where, for each integer $N\geq 2,$,
\begin{align*}
\mathcal{J}^N_{1,i}&:=\{(p,q)\in\Nb^2:i+1\leq p\leq N,\, 1\leq q\leq N+i-p\},\quad i\in\Nb,\\
\mathcal{J}^N_{2,i}&:=\{(p,q)\in\Nb^2:1\leq p\leq N,\, q=i\},\quad i\in\Nb_+,\\
\mathcal{J}^N_{3,i}&:=\{(p,q)\in\Nb^2:1\leq q\leq i-1,\, i-q\leq p\leq N\},\quad i\in\Nb_2,\\
\mathcal{J}^N_{4,i}&:=\{(p,q)\in\Nb^2:p=i,\, 1\leq q\leq N-1\},\quad i\in\Nb_+,
\end{align*}
and, for each integer $N\geq 2$  and $(p,q)\in\mathcal{J}^N_{j,i},$
\begin{align*}
\mathcal{K}^N_{1,i}(p,q)&:=\{p-i\},\quad i\in\Nb,\\
\mathcal{K}^N_{2,i}(p,q)&:=\{k\in\Nb: 1\leq k\leq \min(p,N-i)\},\quad i\in\Nb_+,\\
\mathcal{K}^N_{3,i}(p,q)&:=\{i-q\},\quad i\in\Nb_2,\\
\mathcal{K}^N_{4,i}(p,q)&:=\{k\in\Nb: 1\leq k\leq \min(N-q,i)\},\quad i\in\Nb_+.
\end{align*}
Observe that, for any integer $N\geq 2,$ $\mathcal{J}^N_{j,i}\subsetneqq\mathcal{J}^{N+1}_{j,i}.$ Define 
$\mathcal{J}_{j,i}:=\bigcup_{N=2}^\infty\mathcal{J}^N_{j,i}.$ Since, for each $(p,q)\in\mathcal{J}_{j,i},$ there is $N_0$ such that, for all $N\geq N_0,$
$(p,q)\in\mathcal{J}^N_{j,i}$ we can define $\mathcal{K}_{j,i}(p,q):=\bigcup_{N=N_0}^\infty\mathcal{K}^N_{j,i}(p,q).$ It is not difficult to check that,
$$
\hat{Q}_{j,i}(\beta)=\sum_{(p,q)\in\mathcal{J}_{j,i}}\;\sum_{k\in\mathcal{K}_{j,i}(p,q)}a(p,q;k)\beta_p\beta_q,
$$
for any $\beta\in X_{0,1}.$

\medskip

Let $M,N$ be integers such that $2\leq M\leq N,$ and consider the following inequality,
\begin{equation}\label{QMN}
\begin{split}
|\hat{Q}_{j,i}^N(c^N)-\hat{Q}_{j,i}(c)|_{L^1(0,T)} & \\
& \hspace*{-1cm}\leq \sum_{(p,q)\in\mathcal{J}^M_{j,i}}\;\sum_{k\in\mathcal{K}^M_{j,i}(p,q)}a(p,q;k)\big|c^N_pc^N_q-c_pc_q\big|_{L_1(0,T)}\\
& \hspace*{-1cm}\qquad +\Big|\sum_{(p,q)\in\mathcal{J}^N_{j,i}\setminus\mathcal{J}^M_{j,i}}\;\sum_{k\in\mathcal{K}^N_{j,i}(p,q)}a(p,q;k)c^N_pc^N_q\Big|_{L_1(0,T)}\\
& \hspace*{-1cm}\qquad+\Big|\sum_{(p,q)\in\mathcal{J}^N_{j,i}\setminus\mathcal{J}^M_{j,i}}\;\sum_{k\in\mathcal{K}^N_{j,i}(p,q)}a(p,q;k)c_pc_q\Big|_{L_1(0,T)}
\end{split}
\end{equation}
Since the first sum runs over a fixed finite set of indices, by applying the Lebesgue bounded convergence we can conclude that, for $j\in\{1,2,3,4\}$ and $i\in\Nb_+$, also for $(j,i)=(1,0),$ and excluding $(j,i)=(3,1),$
\begin{equation}\label{lim1}
\lim_{N\to +\infty}\sum_{(p,q)\in\mathcal{J}^M_{j,i}}\;\sum_{k\in\mathcal{K}^M_{j,i}(p,q)}a(p,q;k)\big|c^N_pc^N_q-c_pc_q\big|_{L_1(0,T)}=0.
\end{equation}
Since, in all cases, $\mathcal{K}^N_{j,i}(p,q)\subsetneqq\{k\in\Nb:1\leq k\leq p\},$ we obtain
\begin{align}\label{Sumk}
\sum_{k\in\mathcal{K}^N_{j,i}(p,q)}a(p,q;k)&\;\leq\; \sum_{k=1}^p a(p,q;k)
\;\leq\; C\sum_{k=1}^p (p-k+1)(q+k)q_{p,k}\nonumber\\
&= C\Big(q\sum_{k=1}^p (p-k+1)q_{p,k}+\sum_{k=1}^p (p-k+1)kq_{p,k}\Big)\nonumber \\
&\leq 2C\mathcal{Q}pq.
\end{align}
Also we have, for $\max(2,i)\leq M< N,$
\begin{align*}
\mathcal{J}^N_{2,i}\setminus\mathcal{J}^M_{2,i}&=\big([M+1, N]\cap\Nb\big)\times\{i\},\\
\mathcal{J}^N_{3,i}\setminus\mathcal{J}^M_{3,i}&=\big([M+1,N]\cap\Nb\big)\times\big([1,i-1]\cap\Nb\big),\\
\mathcal{J}^N_{4,i}\setminus\mathcal{J}^M_{4,i}&=\{i\}\times \big([M,N-1]\cap\Nb\big),
\end{align*}
The case $j=1$ has to be tackled in a different way. By using \eqref{Sumk}, we have, for $j\in\{2,3,4\}$ 
\begin{align*}
\Big|\sum_{(p,q)\in\mathcal{J}^N_{j,i}\setminus\mathcal{J}^M_{j,i}}\;\sum_{k\in\mathcal{K}^N_{j,i}(p,q)}a(p,q;k)&c^N_pc^N_q\Big|_{L_1(0,T)}
\\&\leq 2C\mathcal{Q}
\Big|\sum_{(p,q)\in\mathcal{J}^N_{j,i}\setminus\mathcal{J}^M_{j,i}}(pc^N_p)(qc^N_q)\Big|_{L_1(0,T)}\\
&\leq 2C\mathcal{Q}i\|c_0\|\Big|\sum_{p=M}^N pc^N_p\Big|_{L_1(0,T)}\\
&\leq 2C\mathcal{Q}i\|c_0\| \sup_{p\geq M}\frac{p}{\sigma_0(p)}\Big|\sum_{p=M}^N \sigma_0(p)c^N_p\Big|_{L_1(0,T)}
\end{align*}
Therefore, there is a constant $C_0>0$ only depending on $i,T, \mathcal{S}_0$ and $c_0$, such that
\begin{equation}\label{JNpsp}
\Big|\sum_{(p,q)\in\mathcal{J}^N_{j,i}\setminus\mathcal{J}^M_{j,i}}\;\sum_{k\in\mathcal{K}^N_{j,i}(p,q)}a(p,q;k)c^N_pc^N_q\Big|_{L_1(0,T)}\leq C_0\sup_{p\geq M}\frac{p}{\sigma_0(p)}
\end{equation}
The last term in \eqref{QMN} can be dealt in the same way to conclude that
\begin{equation}\label{Jpsp}
\Big|\sum_{(p,q)\in\mathcal{J}_{j,i}\setminus\mathcal{J}^M_{j,i}}\;\sum_{k\in\mathcal{K}_{j,i}(p,q)}a(p,q;k)c_pc_q\Big|_{L_1(0,T)}\leq C_0\sup_{p\geq M}\frac{p}{\sigma_0(p)}.
\end{equation}
By \eqref{QMN}, \eqref{lim1}, \eqref{JNpsp} and \eqref{Jpsp} we have
$$
\limsup_{N\to+\infty}|\hat{Q}_{j,i}^N(c^N)-\hat{Q}_{j,i}(c)|_{L^1(0,T)}\leq 2C_0\sup_{p\geq M}\frac{p}{\sigma_0(p)},
$$
for all $M\geq \max(2,i).$ Since $\sigma_0\in \mathcal{E},$ the right-hand side of the above estimate converges to zero, as $M\to+\infty,$ therefore proving \eqref{claim2} for $j\in\{2,3,4\}$.

For $j=1,$ we consider the following decomposition, for $i\in\Nb$
$$
\mathcal{J}^N_{1,i}\setminus \mathcal{J}^M_{1,i}=J_I\cup J_{II},
$$
where (see Figure~\ref{fig3}),
\begin{align*}
J_I&:=\{(p,q)\in\Nb^2:M+1\leq p\leq N,\; 1\leq q\leq N+i-p\},\\
J_{II}&:=\{(p,q)\in\Nb^2:i+1\leq p\leq M,\; M+1+i-p\leq q\leq N+i-p\}.
\end{align*}

%
%
%
\begin{figure}[h]
	\psfrag{p}{$p$}
	\psfrag{q}{$q$}
	\psfrag{1}{$1$}	
	\psfrag{i}{$i$}
	\psfrag{i+1}{$i+1$}
	\psfrag{M}{$M$}
	\psfrag{M+1}{$M+1$}
	\psfrag{M+i}{$M+i$}	
	\psfrag{p+q=M+1+i}{$p+q=M+1+i$}
	\psfrag{N}{$N$}
	\psfrag{N+i}{$N+i$}
	\psfrag{p+q=N+i}{$p+q=N+i$}
	\psfrag{J_I}{$J_I$}
	\psfrag{J_II}{$J_{II}$}
	\includegraphics[scale=0.40]{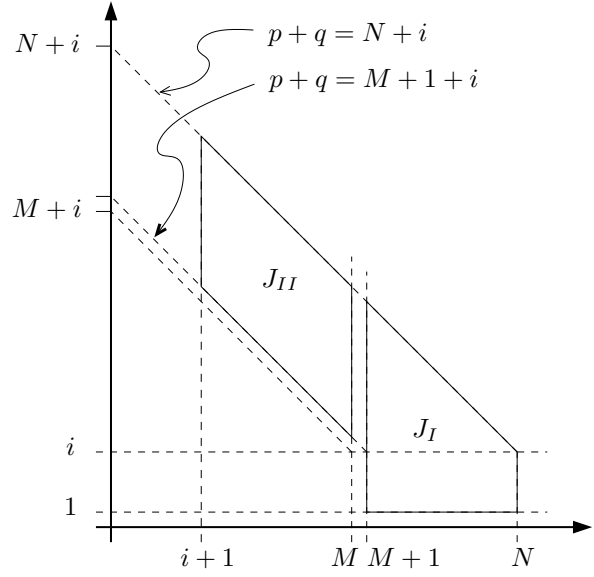}
	\caption{Regions $J_I$ and $J_{II}$ defined in the text.}\label{fig3}
\end{figure}
%
%
%

Then,
\begin{align*}
\Big|\sum_{(p,q)\in J_I}\;\sum_{k\in\mathcal{K}_{1,i}(p,q)}a&(p,q;k)c^N_pc^N_q\Big|_{L_1(0,T)}=
\Big|\sum_{p=M+1}^N\sum_{q=1}^{N+i-p}a(p,q;p-i)c^N_pc^N_q\Big|_{L_1(0,T)}\\
&\leq C(i+1)\Big|\sum_{p=M+1}^N\sum_{q=1}^{N+i-p}(q+p-i)q_{p,p-i}c^N_pc^N_q\Big|_{L_1(0,T)}\\
&\leq 2C\mathcal{Q}\|c_0\|(i+1)\Big|\sum_{p=M+1}^Npc^N_p\Big|_{L_1(0,T)}\\
&\leq 2C\mathcal{Q}\|c_0\|(i+1)\sup_{p\geq M+1}\frac{p}{\sigma_0(p)}\Big|\sum_{p=M+1}^N\sigma_0(p)c^N_p\Big|_{L_1(0,T)},
\end{align*}
so that, there is $C_1>0$ only depending on $i,T,\mathcal{S}_0$ and $c_0$ such that, for all $M\geq i+1,$
\begin{equation}\label{JN1psp}
\Big|\sum_{(p,q)\in J_I}\;\sum_{k\in\mathcal{K}^N_{j,i}(p,q)}a(p,q;k)c^N_pc^N_q\Big|_{L_1(0,T)}\leq C_1\sup_{p\geq M+1}\frac{p}{\sigma_0(p)}.
\end{equation}
Similarly, it is easy to obtain,
\begin{equation}\label{J1psp}
\Big|\sum_{(p,q)\in J_I}\;\sum_{k\in\mathcal{K}^N_{j,i}(p,q)}a(p,q;k)c_pc_q\Big|_{L_1(0,T)}\leq C_1\sup_{p\geq M+1}\frac{p}{\sigma_0(p)}.
\end{equation}
The estimate over $J_{II}$ is more involved:
\begin{align*}
\Big|\sum_{(p,q)\in J_{II}}\;\sum_{k\in\mathcal{K}_{1,i}(p,q)}&a(p,q;k)c^N_pc^N_q\Big|_{L_1(0,T)}=
\Big|\sum_{p=i+1}^M\sum_{q=M+1+i-p}^{N+i-p}a(p,q;p-i)c^N_pc^N_q\Big|_{L_1(0,T)}\\
&\leq C(i+1)\Big|\sum_{p=i+1}^M\sum_{q=M+1+i-p}^{N+i-p}(q+p-i)q_{p,p-i}c^N_pc^N_q\Big|_{L_1(0,T)}\\
&\leq C\mathcal{Q}(i+1)\sup_{r\geq M+1}\frac{r}{\sigma_0(r)}
\Big|\sum_{p=i+1}^M\sum_{q=M+1+i-p}^{N+i-p}\sigma_0(q+p-i)c^N_pc^N_q\Big|_{L_1(0,T)},
\end{align*}
where we have used the fact that for $(p,q)\in J_{II},$ $r:=q+p-i\geq M+1.$

Using the monotonicity of $\sigma_0$ and the inequality \eqref{Convex1} with $m_{\sigma_0}=2$, we get
$$
\sigma_0(q+p-i)\leq \sigma_0(q+p)\leq 3\sigma_0(p)+3\sigma_0(q).
$$
This implies that,
\begin{align*}
\Big|\sum_{p=i+1}^M\sum_{q=M+1+i-p}^{N+i-p}\sigma_0(q+p-i)c^N_pc^N_q\Big|_{L_1(0,T)}
&\leq 3\Big|\sum_{p=i+1}^M\sum_{q=M+1+i-p}^{N+i-p}(\sigma_0(p)+\sigma_0(q))c^N_pc^N_q\Big|_{L_1(0,T)}\nonumber\\
&\leq 6\mathcal{S}_0\|c_0\|T.
\end{align*}
Hence, there is a constant $C_2>0$ only depending on $i,T,\mathcal{S}_0$ and $\|c_0\|$, such that
\begin{equation}\label{JN2psp}
\Big|\sum_{(p,q)\in J_{II}}\;\sum_{k\in\mathcal{K}^N_{j,i}(p,q)}a(p,q;k)c^N_pc^N_q\Big|_{L_1(0,T)}\leq C_2\sup_{r\geq M+1}\frac{r}{\sigma_0(r)},
\end{equation} 
and, obtained in the same way,
\begin{equation}\label{J2psp}
\Big|\sum_{(p,q)\in J_{II}}\;\sum_{k\in\mathcal{K}_{j,i}(p,q)}a(p,q;k)c_pc_q\Big|_{L_1(0,T)}\leq C_2\sup_{r\geq M+1}\frac{r}{\sigma_0(r)}.
\end{equation} 
As before, by \eqref{lim1}, \eqref{JN1psp}, \eqref{J1psp}, \eqref{JN2psp} and \eqref{J2psp}, we obtain, for $i\in\Nb,$
$$
\limsup_{N\to+\infty}|\hat{Q}_{1,i}^N(c^N)-\hat{Q}_{1,i}(c)|_{L^1(0,T)}\leq C_3\sup_{p\geq M+1}\frac{p}{\sigma_0(p)},
$$
for some constant $C_3>0$ only depending on $i, T, \mathcal{S}_0$ and $c_0.$ Since $\sigma_0\in \mathcal{E},$ the right-hand side of the above estimate converges to zero, as $M\to+\infty,$ therefore proving \eqref{claim2} for $j=1$.

\medskip

Now we tackle the terms whose rate coefficients are of fragmentation-type, namely: $a(p,0;k).$ 
Define $A^N_{j,i}$ as the sum of terms in $Q^N_{j,i}$ which are of this type. According to \eqref{hatQN1I}-\eqref{hatQN4I} we have, for $N\geq i+2,$
\begin{align}\label{AN1I}
A^N_{1,i}(c^N)&=\displaystyle (1-\delta_{i,0})\sum_{k=1}^{N-i}a(i+k,0;k)c^N_{i+k}c^N_0,\qquad\;\;\; i\in\{0, 1,\dots,N-1\},\\
\label{AN2I}A^N_{2,i}(c^N)&=\displaystyle -\delta_{i,0}\sum_{k=1}^{N-1}\sum_{j=k+1}^Na(j,0;k)c^N_jc^N_0,\qquad\qquad\; i\in\{0, 1,\dots,N-1\},\\
\label{AN3I}A^N_{3,i}(c^N)&=\sum_{j=i+1}^Na(j,0,i)c^N_jc^N_0,\quad\qquad\qquad\;\;\;\, i\in\{1,\dots,N\}, \\
\label{AN4I}A^N_{4,i}(c^N)&=-\sum_{k=1}^{i-1}a(i,0;k)c^N_0c^N_i\,,\qquad\qquad\;\; \quad i\in\{1,\dots,N\},
\end{align}
where $\delta_{i,0}$ is the Kronecker symbol.
Here, we used Definition~\eqref{nullcoef} that all coefficients of type $a(k,0; k)$ are zero. 
Correspondingly, we define $A_{j,i}$ as the sums of terms in $Q_{j,i}$ with rate coefficients of the same type which are obtained from 
\eqref{AN1I}-\eqref{AN4I} making $N\to\infty$ in the upper limits of the sums. Our next claim is that $A_{j,i}(c)\in L_1(0,T)$ and 
\begin{equation}\label{claim3}
\lim_{N\to +\infty}|A_{j,i}^N(c^N)-A_{j,i}(c)|_{L^1(0,T)}=0.
\end{equation}
Observe that, for each $i=1,\dots,N,$ the sum in \eqref{AN4I} runs over a finite set of indices and that
$$
\int_0^T|A^N_{4,i}(c^N(t))-A_{4,i}(c(t))|dt=\alpha(i)\int_0^T|c_0^N(t)c_i^N(t)-c_0(t)c_i(t)|dt,
$$  
and hence, our claim is easily proved for $j=4$ by the bounded convergence theorem. For the other cases we write, for $j=1,2,3,$ 
$$
A^N_{j,i}(c)=c^N_0\sum_{(p,k)\in \mathcal{I}_{j,i}^N}a(p,0;k)c_p^N,
$$
where,
\begin{align*}
\mathcal{I}_{1,i}^N&:=\{(p,k)\in\Nb^2: i+1\leq p\leq N,\; k=p-i \},\quad\, i\in\{1,\dots,N-1\},\\
\mathcal{I}_{2,0}^N&:=\{(p,k)\in\Nb^2: 2\leq p\leq N,\; 1\leq k\leq p-1\},\\
\mathcal{I}_{3,i}^N&:=\{(p,k)\in\Nb^2: i+1\leq p\leq N,\; k=i \},\quad\qquad i\in\{1,\dots,N\}.
\end{align*}
Observe that $\mathcal{I}_{j,i}^N\subsetneqq \mathcal{I}_{j,i}^{N+1}.$ Also, define $\mathcal{I}_{j,i}:=\bigcup_{N=i+1}^\infty\mathcal{I}_{j,i}^N,$ and consider the sets 
 $$
 \mathcal{J}^N:=\{(p,k)\in(\Nb_+)^2:1\leq p\leq N,\; 1\leq k\leq p-1\},
 $$
and $\mathcal{J}:=\bigcup_{N=1}^\infty\mathcal{J}^N.$ Let us also recall the sets $\mathcal{J}^N_0$  from Lemma \ref{HMomentLem1} and consider $\mathcal{J}_0:=\bigcup_{N=M_0+1}^\infty\mathcal{J}_0^N$. Then,  the set $\mathcal{J}\setminus\mathcal{J}_0$ is finite. Therefore,
\begin{equation}\label{sum1L1}
\sum_{(p,k)\in\mathcal{J}\setminus\mathcal{J}_0} a(p,0;k)c_0c_p\in L_1(0,T).
\end{equation}
On the other hand, for $(p,k)\in\mathcal{J}_0^N$ we have, by Lemma \ref{HMomentLem1}, and the fact that, for $p\geq M_0+1,$ $\sigma_0(p-1)/p-1\geq 1,$ for fixed $\mu\geq M_0+1,$
$$
\int_0^T\sum_{(p,k)\in\mathcal{J}_0^\mu} a(p,0;k)c^N_0(t)c^N_p(t)dt\leq 
\gamma_T \sum_{i=0}^\infty \sigma_0(i) c_{0 i}.
$$
By letting $N\to +\infty,$ by using the bounded convergence theorem, we can replace in the previous inequality $c^N_p$ by $c_p.$ Then, by the monotone convergence theorem, by letting $\mu\to +\infty$, we obtain 
\begin{equation}\label{sum2L1}
\sum_{(p,k)\in\mathcal{J}_0} a(p,0;k)c_0c_p\in L_1(0,T).
\end{equation}
By \eqref{sum1L1} and \eqref{sum2L1}, we conclude that,
$$
\sum_{(p,k)\in\mathcal{J}} a(p,0;k)c_0c_p\in L_1(0,T)
$$
Since, for $j=1,2,3,$ and corresponding indices $i,$  $\mathcal{I}_{j,i}\subset\mathcal{J},$ we can conclude that,
\begin{equation}
A_{j,i}(c)\in L_1(0,T).
\end{equation}
This, together with \eqref{claim1} allows us to conclude that
$$
Q_{j,i}(c)\in L_1(0,T),\quad j\in\{1,2,3,4\}\quad (\text{or } j\in\{1,2\}\text{ if }i=0),
$$
which is condition (ii) of Definition \ref{Def Sol}.

The proof of \eqref{claim3} proceeds in a way similar to the proof of claim \eqref{claim2}, by considering the inequality corresponding to \eqref{QMN} but now with the $\hat{Q}$'s replaced by the $A$'s and with the obvious replacements on the right-hand side. So, for each pair of integers $M,N$ such that $i+1\leq M\leq N,$ we have
\begin{align*}
\mathcal{I}_{1,i}^N\setminus\mathcal{I}_{1,i}^M&=\{(p,k)\in\Nb^2: M+1\leq p\leq N,\; k=p-i \},\quad i\in\{1,\dots,N-1\},\\
\mathcal{I}_{2,0}^N\setminus\mathcal{I}_{2,0}^M&=\{(p,k)\in\Nb^2: M+1\leq p\leq N,\; 1\leq k\leq p-1\},\\
\mathcal{I}_{3,i}^N\setminus\mathcal{I}_{3,i}^M&=\{(p,k)\in\Nb^2: M+1\leq p\leq N,\; k=i \},\quad i\in\{1,\dots,N\}.
\end{align*}
By the bounded convergence theorem, for each choice of $M$ we have,
\begin{equation}\label{lim2}
\lim_{N\to +\infty}\sum_{(p,k)\in\mathcal{I}^M_{j,i}}a(p,0;k)\big|c^N_pc^N_0-c_pc_0\big|_{L_1(0,T)}=0.
\end{equation}
With reference to Lemma \ref{HMomentLem1}, we have, for $j=1,2,3$ and $i$ in the corresponding sets, if $N\geq M\geq M_0,$ then $\mathcal{I}^N_{j,i}\setminus\mathcal{I}^M_{j,i}\subsetneqq\mathcal{J}_0$ so that \eqref{gamma1lemma1} is true for $\mathcal{J}_0$ replaced by $\mathcal{I}^N_{j,i}\setminus\mathcal{I}^M_{j,i}$, that is, there is a constant $C(T)>0,$ such that,
\begin{equation}\label{ineqji}
\int_0^T \sum_{(p,k)\in\mathcal{I}^N_{j,i}\setminus\mathcal{I}^M_{j,i}}\frac{\sigma_0(p-1)}{p-1}a(p,0;k)c^N_0(t)c^N_p(t)dt
\leq C(T).
\end{equation}

For $j=1,2,3$ and $(p,k)\in\mathcal{I}^N_{j,i}\setminus\mathcal{I}^M_{j,i}$, we have $p\geq M+1$ and therefore, given $\varepsilon\in(0,1),$ for $M$ sufficiently large we have
$$
\sup_{(p,k)\in\mathcal{I}^N_{j,i}\setminus\mathcal{I}^M_{j,i}}\frac{\sigma_0(p-1)}{p-1}= \frac{\sigma_0(M)}{M}\geq \frac{2C(T)}{\varepsilon},
$$
which, by \eqref{ineqji} allows us to conclude that,
\begin{equation}\label{ineqIN2}
\Big|\sum_{(p,k)\in\mathcal{I}^N_{j,i}\setminus\mathcal{I}^M_{j,i}}a(p,0;k)c^N_0c^N_p\Big|_{L_1(0,T)}
\leq \frac{\varepsilon}{2} .
\end{equation}
Proceeding similarly, we obtain
\begin{equation}\label{ineqI2}
\Big|\sum_{(p,k)\in\mathcal{I}_{j,i}\setminus\mathcal{I}^M_{j,i}}a(p,0;k)c_0c_p\Big|_{L_1(0,T)}
\leq \frac{\varepsilon}{2} .
\end{equation}
From \eqref{lim2}, \eqref{ineqIN2} and \eqref{ineqI2}, we obtain
$$
\limsup_{N\to +\infty}|A_{j,i}^N(c^N)-A_{j,i}(c)|_{L^1(0,T)}\leq \varepsilon,
$$
from which, by letting $\varepsilon\to 0,$ we obtain \eqref{claim3}.

From \eqref{claim2} and \eqref{claim3} we prove that 
 for $j\in\{1,2,3,4\},$ $i\in\Nb_+,$ and also for $(j,i)=(1,0),$
\begin{equation}\label{QNconvQ}
\lim_{N\to +\infty}|Q_{j,i}^N(c^N)-Q_{j,i}(c)|_{L^1(0,T)}=0.
\end{equation}
Then, from the integrated version of the truncated system \eqref{NDGED}-\eqref{Q4Ni}, and also Corollary~\ref{Corollary1}, \eqref{Converge1} and \eqref{normest}, together with \eqref{QNconvQ}, it is straightforward to obtain that $c$ satisfies condition (iii) of Definition~\ref{Def Sol}.
The continuity of $c(\cdot)$ is a consequence of this.       \end{Proofof}
\medskip

Analogously with what is sometimes the practice in the coagulation-fragmenta\-tion literature \cite{carr92}
we shall call \emph{admissible} any solution to the initial value problem \eqref{DGED}, \eqref{init} obtained 
as limit of solutions to the \eqref{NDGED}--\eqref{Tinit} when $N\to\infty$, in the sense used in the above proof 
of Theorem~\ref{ExistThm1}.

\medskip

The proof that, in the isolated case, admissible solutions of \eqref{DGED} conserve the moments \eqref{moments} with
$r=0$ (total number of clusters) and $r=1$ (total mass of clusters) is easily done  using the fact
that if an initial condition $c_0$ is in $X^+_{0,1}$ then it satisfies \eqref{sigmac0}, for
some $\sigma_0\in\mathcal{E}$, and so it is slightly more regular,
and this extra regularity is inherited by the admissible solution (remember 
Lemma~\ref{HMomentLem1} and \eqref{sigmact}).


\begin{Proofof} of theorem \ref{TPTheorem}: Let $T>0$ and $2 \leq L < N$. Using Corollary~\ref{Corollary1} we write, for all $t\in [0, T).$
\begin{align}
\mathcal{P}_0(t)  - \mathcal{P}_0(0)    
& = \sum_{i=0}^{N} ( c_i(t)   - c_{0 i} )   + \sum_{i=N+1}^\infty c_i(t) - \sum_{i=N+1}^\infty c_{0 i}   \nonumber\\
& = \sum_{i=0}^{N} ( c_i(t)   - c_{i}^N(t) )   +  \sum_{i=0}^{N} ( c_i^N(t)   - c_{0 i} )  + \sum_{i=N+1}^\infty c_i(t) - \sum_{i=N+1}^\infty c_{0 i}   \nonumber\\
& = \sum_{i=0}^{L} ( c_i(t)   - c_{i}^N(t) )   +   \sum_{i=L+1}^{\infty} c_i(t)   - \sum_{i=L+1}^{N}c_{i}^N(t)   - \sum_{i=N+1}^\infty c_{0 i}   \nonumber
\end{align} 
Thus, from the fact that $c_0\in X_{0, 1}$, 
\begin{align}
\abs{\mathcal{P}_0(t)  - \mathcal{P}_0(0)  } 
& \leq  \sum_{i=0}^{L} \left| c_i(t)   - c_{i}^N(t) \right|   +   \sum_{i=L+1}^{\infty} |c_i(t)| +   
\sum_{i=L+1}^{N}|c_{i}^N(t)|   + \sum_{i=N+1}^\infty |c_{0 i}|,    \nonumber   \\
& \leq   \sum_{i=0}^{L} \left| c_i(t)   - c_{i}^N(t) \right| +   \sum_{i=L+1}^{\infty} |c_i(t)| \;\,+ \nonumber \\
& \qquad\qquad\quad +  \frac{1}{L+1}\sum_{i=L+1}^{N}|ic_{i}^N(t)|+ \sum_{i=N+1}^\infty |c_{0 i}|.  \nonumber
\end{align}
Now letting $N\to \infty$ and using \eqref{Converge1} and \eqref{sigmactlim} it follows that 
\begin{align}
\abs{\mathcal{P}_0(t)  - \mathcal{P}_0(0)  } 
& \leq     \sum_{i=L+1}^{\infty} |c_i(t)| + \frac{1}{L+1} \gamma_T\mathcal{S}_0,  
\end{align}
and letting $L\to \infty$ we conclude that $\abs{\mathcal{P}_0(t)  - \mathcal{P}_0(0)  } = 0,$ thus proving 
\eqref{P0conserve}.

\medskip

To prove \eqref{P1conserve} repeat the computations above, now for $\mathcal{P}_1(t)  - \mathcal{P}_1(0)$.
We get
\begin{align}
\abs{\mathcal{P}_1(t)  - \mathcal{P}_1(0)  } 
& \leq  \sum_{i=0}^{L} i\left| c_i(t)   - c_{i}^N(t) \right|   +   \sum_{i=L+1}^{\infty} |ic_i(t)| +   
\sum_{i=L+1}^{N}|ic_{i}^N(t)|   + \sum_{i=N+1}^\infty |ic_{0 i}|,    \nonumber   \\
& \leq   \sum_{i=0}^{L} i\left| c_i(t)   - c_{i}^N(t) \right| +   \sum_{i=L+1}^{\infty} |ic_i(t)| \;\,+ \nonumber \\
& \qquad\qquad\quad +   \frac{L+1}{\sigma_0(L+1)}\sum_{i=L+1}^{N}|\sigma_0(i)c_{i}^N(t)|+ \sum_{i=N+1}^\infty |ic_{0 i}|   \nonumber
\end{align}
and now repeating the process of letting first $N\to\infty$, then $L\to\infty$, and using \eqref{sigmactlim}
and \eqref{Convex}, we conclude \eqref{P1conserve}.\end{Proofof}

\subsection{The non-isolated case}
In the non-isolated case theorems~\ref{ExistThm1} and \ref{TPTheorem}
can be stated almost exactly verbatim (the only difference 
being the lack of conservation of the cluster number density)
with essentially the same proofs, where the only difference is that the terms $Q_{j,0}$ are now all zero. For completeness we state them now, omitting the proofs.


\section{A regularity result}\label{secreg}

All along this work we have been considering conditions \eqref{qikbound}-\eqref{bound} with the only condition imposed on the rate coefficients of fragmentation-type, $a(p,0;k),$ being their nonnegativity, besides the ``physical'' conditions \eqref{nullcoef}-\eqref{symfrag}.
This lack of upper bounds prevented us from using Ascoli-Arzel\`a theorem in the proof of the existence theorems. 
It also prevent us to obtain uniform convergence properties for the series defining each $Q_{j,i}(c(\cdot)),$ for the solution $c$ obtained in the existence theorem and therefore, to obtain more regularity for this solution than its continuity. However, if we extend the upper bound \eqref{bound} to $j=0,$ we obtain the following regularity result which, in particular, proves that the solution constructed in Theorem~\ref{ExistThm1} is a solution of \eqref{DGED}-\eqref{init} in the classical sense:


\begin{Proofof} of the theorem \ref{Regthm}:
Similar to the function $\Phi$ introduced in the proof of Theorem~\ref{ExistThm1}, we define the function $\Psi$ by
$$
\Psi(\beta):=\sum_{(i,j,k)\in\mathcal{I}_0}a(i,j;k)\beta_i\beta_j,
$$
where,
$$
\mathcal{I}_0:=\{(i,j,k)\in\Nb^3 : 1\leq k\leq i\},
$$
the domain of which, $D_\Psi,$ is the subset of $X^+_{0,1}$ formed by the elements $\beta$ for which the above series is convergent. By proceeding as in the proof of Theorem \ref{ExistThm1}, we obtain, for any integer $\eta\geq 2,$
\begin{equation}\label{Psibeta}
\sum_{i=1}^\eta\sum_{j=0}^\eta\sum_{k=1}^i a(i,j;k)\beta_i\beta_j\leq 2C\mathcal{Q}
\Big(\sum_{i=1}^\eta i\beta_i\Big)\Big(\sum_{i=1}^\eta j\beta_j + \sum_{j=0}^\eta \beta_j\Big).
\end{equation}
This shows that, in fact, $D_\Psi=X^+_{0,1}.$
On the other hand, from the assumption that $\sum_{i=1}^\infty ic_i(t)$ is constant for $t\in[0,T)$ Dini's theorem 
allow us to conclude that this series of nonnegative functions is uniformly convergent on any compact subset of $[0,T)$. And so $\sum_{i=0}^\infty c_i(t)$ is also uniformly convergent on compact subsets of $[0, T)$. 
This, together with \eqref{Psibeta} ensure that the same is true for the series defining $\Psi(c(\cdot))$.  
Since, for $j=1,2,3,4$ with $i\in\mathbb{N}_+$, and $j=1,2$ with $i=0$, the the indices of sums defining 
$Q_{j,i}(c)$ run over subsets of $\mathcal{I}_0,$ we conclude that these are also uniformly convergent on any compact subset of $[0,T)$ and therefore, by the continuity of $c,$ we conclude that
$Q_{j,i}(c(\cdot))$ is continuous in $[0,T)$. But then, by (iii) in Definition~\ref{Def Sol}, the 
thesis of the lemma is obtained.
\end{Proofof}
\begin{remark}
We remark that the proof presented above is valid for
both the isolated and non-isolated cases as it is based on the
bound \eqref{Psibeta} and the conservation of mass (together
with Dini's theorem), both of which hold true in both regimes.
\end{remark}
\begin{remark}
In the non-isolated case discussed in sections \ref{GEDGM} and \ref{MS}, when the DGED system is reduced to the standard coagulation-fragmentation equations, if we impose the conditions of Theorem~\ref{Regthm} we have for the fragmentation coefficients,
\begin{align*}
b_{j,k}&=2a(j+k,0;k)c_0\\
&\leq 2C(j+1)kq_{j+k,k}c_0\\
&\leq Kjk,
\end{align*}
with $K=2(C+1)\mathcal{Q}c_0,$ and we recover Theorem~5.2 of \cite{Ball:1990}.
\end{remark}

\section{Uniqueness of solutions}\label{uniq}

In this section, 
we 
prove 
the partial 
uniqueness result for \eqref{DGED}--\eqref{init} 
stated in theorem \ref{Uniqthm}.
The proof is based on ideas commonly used in studies of 
coagulation-type equations (see, e.g., \cite{Ball:1990,PL:2002}) and so we will skip most of the details
of the computations. Those ideas are the following: we assume that, for some $T\in (0, \infty)$ the
initial value problem \eqref{DGED}--\eqref{init} has two solutions $c=(c_i)$ and $d=(d_i)$.
Defining $x:=c-d$, we prove that, on $[0, T]$,  this function satisfies a differential inequality
\[
\psi(t) \leq \int_0^t \psi(s)ds,
\]
where $\psi(t) := \sum_{i=0}^\infty (1+i^\alpha)|x_i(t)|.$ From this, by Gronwall's inequality, we get
$x_i(t)\equiv 0$, and hence uniqueness follows.

\medskip

To implement this idea we need to get an evolution equation for solutions to \eqref{DGED} similar
to \eqref{gNmom}, in Proposition- \ref{Prop1}, which was valid for solutions of the truncated system 
\eqref{NDGED}--\eqref{Q4Ni}. This is stated in the next proposition, whose proof is similar to that of
Proposition~\ref{Prop1} and we leave it to the reader.

\begin{prop}\label{Prop1general}
Let $c=(c_i)_{i\geq 0 }$ be any solution of \eqref{DGED} in the isolated case. 
Then, for every $n\in\Nb$ and sequence $(g_i)$ we have
\begin{align}
\sum_{i=0}^ng_ic_i(t) & = \sum_{i=0}^ng_ic_i(0)  \nonumber   \\
& \quad + \int_0^t\sum_{k=1}^{n}\sum_{i=k}^n\sum_{j=0}^{n-k}(g_{j+k}+g_{i-k}-g_j-g_i)a(i, j; k)c_i(s)c_j(s)ds   \label{gmom1}  \\
& \quad + \int_0^t\sum_{j=1}^4 R_{j,n}(c(s))ds \label{gmom2}
\end{align}
where
\begin{align}
R_{1,n}(c(\cdot)) & := \sum_{k=1}^{n}\sum_{i=n+1}^\infty \sum_{j=0}^{n-k} g_{j+k}a(i, j; k)c_{i}c_{j},		\label{R1N} \\
R_{2,n}(c(\cdot)) & := \Bigl(\sum_{k=1}^\infty \sum_{i=k}^{n+k}\sum_{j=0}^\infty - \sum_{k=1}^n\sum_{i=k}^n\sum_{j=0}^{n-k} \Bigr) g_{i-k} a(i, j; k)c_{i}c_{j},		\label{R2N} \\
R_{3,n}(c(\cdot)) & := - \Bigl(\sum_{k=1}^\infty \sum_{i=k}^{\infty}\sum_{j=0}^n - \sum_{k=1}^n\sum_{i=k}^n\sum_{j=0}^{n-k} \Bigr) g_j a(j, i-k; k)c_{j}c_{i-k},		\label{R3N} \\
R_{4,n}(c(\cdot)) & := - \sum_{k=1}^{n}\sum_{i=k}^n \sum_{j=n-k+1}^{\infty} g_{i}a(i, j; k)c_{i}c_{j}.		\label{R4N} 
\end{align}
\end{prop}
\medskip

Observe that in \eqref{gmom1}, when $k=n$ the only term of the triple sum has $(i,j,k)=(n,0,n)$ and so, either
from $g_{j+k}+g_{i-k}-g_j-g_i = g_n+g_0-g_0-g_n = 0$, or by \eqref{nullcoef}, its contribution to \eqref{gmom1}
is identically zero. This means that the sum in $k$ can be done just from $1$ to $n-1$, exactly as in \eqref{gNmom}.

\medskip

We will assume that the rate coefficients $a(i,j; k)$ satisfy the following: there exist a constant 
$\alpha \in [0, \frac{1}{2})$ and a positive constant
$C$ such that, for all integers $i, k$ such
that $1 \leq k \leq i$, and all $j \in \Nb$, we have
\begin{equation}\label{boundalpha}
a(i, j; k) \leq C (i-k+1)^\alpha(j^\alpha + k^\alpha)q_{i,k},
\end{equation}
and $q_{i,k}$ satisfy the condition \eqref{qikbound} used in the existence result.

\medskip

Note that condition \eqref{boundalpha} is more restrictive than \eqref{bound} not only due to the
constant $\alpha$ being smaller than $1$, but  because it is assumed to be valid also for the case $j=0$
which we had not assumed to be the case in \eqref{bound} (remember the discussion in Section~\ref{secreg}.)

\medskip

We also observe that from \eqref{qikbound} and \eqref{boundalpha} it follows that 
\begin{equation}\label{qikboundalpha}
\sum_{k=1}^i k^\alpha(i-k+1)^\alpha q_{i,k}\leq \mathcal{Q}i^\alpha.
\end{equation}
 In fact, from $\min_{1\leq k\leq i}\{k(i-k+1)\} = i,$ we have
 \begin{align*}
 \mathcal{Q}i 
 & \geq \sum_{k=1}^i k(i-k+1) q_{i,k} 
 = \sum_{k=1}^i \bigl(k(i-k+1)\bigr)^{1-\alpha}k^\alpha(i-k+1)^\alpha q_{i,k} \\
& \geq i^{1-\alpha}\sum_{k=1}^i k^\alpha(i-k+1)^\alpha q_{i,k},
 \end{align*}
 from which \eqref{qikboundalpha} readily follows. 
 
 \medskip
 
We can now prove the uniqueness result stated in theorem \ref{Uniqthm}.
 

\begin{Proofof} of the theorem \ref{Uniqthm}:
Suppose $c=(c_i)$ and $d=(d_i)$ are two solutions of \eqref{DGED}--\eqref{init} on $[0,T].$ Let $x:= c-d$.
Note that, since that $c(0)$ and $d(0)$ are both equal to the initial condition $c_0$, and thus $x(0)=0$. 
Also, for all $i$ and $j$ it is clear that $c_ic_j - d_id_j = c_ix_j + d_jx_i.$
Observing that, for every absolutely continuous function $u(t)$, the function $t\mapsto |u(t)|$ is also
absolutely continuous and $\frac{d}{dt}|u(t)| = \sgn (u(t))\frac{du}{dt}$ almost everywhere 
(where $\sgn(\cdot)$ denotes the sign function),
we can use Proposition~\ref{Prop1general} to write
\begin{align}
\sum_{i=0}^n(1+i^\alpha)|x_i(t)| 
& = \int_0^t\sum_{k=1}^{n}\sum_{i=k}^n\sum_{j=0}^{n-k}\tilde{g}(i, j; k)a(i, j; k)\bigl(c_i(s)x_j(s) + d_j(s)x_i(s)\bigr)ds \label{Un}  \\
& \qquad +  \int_0^t\sum_{j=1}^4 \Delta R_{j,n}(s)ds,  \label{Vn}
\end{align}
where $\tilde{g}(i, j; k) := g_{j+k}+g_{i-k}-g_j-g_i$, with $g_i = (1+i^\alpha)\sgn(x_i)$, and $\Delta R_{j,n}(\cdot) := 
R_{j,n}(c(\cdot)) - R_{j,n}(d(\cdot)).$ 
Writing $x_\ell = \sgn(x_\ell)|x_\ell|$ and remembering that the sign function is either $-1, 0$ or $1$, we can write
\begin{align}
\tilde{g}(i, j; k)x_j & = (g_{j+k}+g_{i-k}-g_j-g_i)x_j \nonumber \\
& \leq \bigl(2 + (j+k)^\alpha+(i-k)^\alpha- j^\alpha +i^\alpha\bigr)|x_j|,  \nonumber
\end{align}
now using the concavity of $u\mapsto u^\alpha$ and the fact that for $0\leq u\leq i$ the function
$u\mapsto u^\alpha + (i-u)^\alpha$, has a maximum at $u= i/2$, whose value is $2^{1-\alpha} i^\alpha$ we get
\begin{align}
\tilde{g}(i, j; k)x_j  & \leq \bigl(2 + j^\alpha+k^\alpha+(i-k)^\alpha- j^\alpha +i^\alpha\bigr)|x_j| \nonumber \\
& \leq \bigl(2 + (2^{1-\alpha}+1\bigr)i^\alpha)|x_j| \nonumber \\
& \leq (2 + 3i^\alpha)|x_j|.     \label{gxj}
\end{align}
Analogously,
\begin{align}
\tilde{g}(i, j; k)x_i & \leq \bigl(2 + (j+k)^\alpha+(i-k)^\alpha+ j^\alpha -i^\alpha\bigr)|x_i|, \nonumber \\
& \leq \bigl(2 + j^\alpha+k^\alpha+(i-k)^\alpha+ j^\alpha -i^\alpha\bigr)|x_i|\nonumber \\
& \leq (2 + 2j^\alpha + k^\alpha)|x_i|. \label{gxi}
\end{align}
Plugging  \eqref{boundalpha}, \eqref{qikboundalpha}, \eqref{gxj} and \eqref{gxi} into \eqref{Un} we can 
bound from above the function in the integral in \eqref{Un}  by
\begin{equation}
7 CQ\sup_{t\in [0, T]}\bigl(\|c(t)\| + \|d(t)\|\bigr)\sum_{i=0}^n(1+i^\alpha)|x_i(s)|. \label{limUn}
\end{equation}

\medskip

To estimate the terms $\Delta R_{j,n}$ we use the same tools as above (assumptions $\alpha \in [0, \frac{1}{2})$, \eqref{boundalpha}, \eqref{qikboundalpha}, the concavity of $i\mapsto i^\alpha$, values of the sign function) and
recalling that $|x_\ell| = |c_\ell-d_\ell| \leq c_\ell + d_\ell,$ to obtain
\begin{align}
\Delta R_{1,n} & \leq \frac{6 CQ}{(n+1)^{1-2\alpha}}\sup_{t\in [0, T]}\bigl(\|c(t)\|+\|d(t)\|\bigr)^2,\label{boundR1n} \\
\Delta R_{2,n} & \leq \frac{12 CQ}{(n+1)^{1-2\alpha}}\sup_{t\in [0, T]}\bigl(\|c(t)\|+\|d(t)\|\bigr)^2 + \nonumber \\
& \qquad 
+ 4CQ\sum_{i=1}^n\sum_{j=n-i+1}^\infty i^{2\alpha}(c_i+d_i)j^{2\alpha}(c_j+d_j) \label{boundR2n}\\
\Delta R_{3,n} & \leq \frac{4CQ}{(n+1)^{1-\alpha}}\sup_{t\in [0, T]}\bigl(\|c(t)\|+\|d(t)\|\bigr)^2 + \nonumber \\
& \qquad 
+ 8CQ\sum_{i=1}^n\sum_{j=n-i+1}^n i^{2\alpha}(c_i+d_i) j^{2\alpha}(c_j+d_j), \label{boundR3n} \\
\Delta R_{4,n} & \leq  4CQ\sum_{i=1}^n\sum_{j=n-i+1}^\infty i^{2\alpha}(c_i+d_i)(c_j+d_j). \label{boundR4n} 
\end{align}
For $s\in [0, T]$ we have $\sum_{i=0}^n(1+i^\alpha)|x_i(s)| \leq \|x(s)\| \leq \sup_{s\in [0, T]}\|x(s)\| < \infty$,
and so the bound \eqref{limUn} and the monotone convergence theorem imply
that the integral \eqref{Un} is bounded above by
\begin{equation*}
7 CQ\sup_{t\in [0, T]}\bigl(\|c(t)\| + \|d(t)\|\bigr)\int_0^t\sum_{i=0}^\infty(1+i^\alpha)|x_i(s)|ds, \label{limlimUn}
\end{equation*}
when $n\to \infty.$

\medskip

In the same way, from \eqref{boundR1n}--\eqref{boundR4n} and the assumption $\alpha < \frac{1}{2}$ we conclude, using the dominated convergence theorem, that  the integral \eqref{Vn} converges to zero as $n\to \infty$.

\medskip

Hence, we can now pass to the limit $n\to \infty$ in \eqref{Un}--\eqref{Vn} obtaining, for $t\in [0, T]$ the 
inequality
\[
\sum_{i=0}^\infty(1+i^\alpha)|x_i(t)| \leq 7 CQ\sup_{t\in [0, T]}\bigl(\|c(t)\| + \|d(t)\|\bigr)\int_0^t\sum_{i=0}^\infty(1+i^\alpha)|x_i(s)|ds,
\]
and thus Gronwall's lemma implies that
\[
\sum_{i=0}^\infty(1+i^\alpha)|x_i(t)| = 0, \forall t\in[0, T],
\]
whence $x_i(t)\equiv 0$, and the uniqueness is proved.
\end{Proofof}
 The reader may have noted that the above proof is valid for both the isolated and the non-isolated cases.

 \medskip

\section{Final Remarks}\label{finalrem}

In this paper we introduced a new model for exchange-driven growth of clusters made up of a discrete number of
particles. This model allows for the exchange of a number $k$ of particles between two clusters thus generalizing
the exchange-driven model existing in the literature, for which only a single particle ($k=1$) can be exchanged
in each reaction between clusters.

\medskip

The mathematical description of the time evolution of the concentration at time $t$ of 
clusters made of $i\geq 0$ identical particles, $c_i(t),$ is the system of ordinary differential equations 
\eqref{DGED}, which we called the Discrete Generalized Exchange Driven system (DGED). In this paper, 
after introducing the model and relating it with other cluster models (exchange-driven, coagulation-fragmentation)
we proved an existence result under reasonably general conditions on the rate coefficients, and two
conservation results. Under a slightly more 
restrictive hypothesis on the rate coefficients, we prove uniqueness of solutions to the initial value problems. 
Also under more restrictive assumptions than those used for the existence proof we prove a regularity result.

\medskip

One aspect that we did not touch in this study is the long-time behaviour of solutions. This we expect to be able
to tackle in future work. In this direction, proving that the set of solutions constitute a semi-group in an appropriate
topology will be crucial. Furthermore, existing studies about the exchange-driven and the coagulation-fragmentation systems
suggest that under appropriate conditions DGED can have nonzero equilibria and the dynamic convergence of solutions to these
equilibria as $t\to\infty$ can be studied using a Lyapunov function. We conclude this paper with some further
comments about this aspect.

\medskip

System \eqref{DGED} can be rewritten in a way that highlights the balance of reversible chemical reactions between clusters;
namely, by grouping together $Q_{1,i}(c)+Q_{2,i}(c)$ we model reactions of the type $\langle i+k\rangle + \langle j\rangle\rightleftharpoons \langle i\rangle + \langle j+k\rangle$, and looking at $Q_{3,i}(c)+Q_{4,i}(c)$ 
we are considering the reactions $\langle i-k\rangle + \langle j\rangle\rightleftharpoons \langle i\rangle + \langle j-k\rangle.$
Thus, we can write \eqref{DGED} as follows
\begin{equation}
\dot{c}_i = \sum_{k=1}^i W_{i-k;k}(c) - \sum_{k=1}^\infty W_{i;k}(c)\label{DGEDbalanced}
\end{equation} 
with 
\begin{equation}
W_{p;k}(c):= \sum_{q=k}^\infty \omega(q,p;k)(c),\label{macrobalance}
\end{equation}
and 
\begin{equation}
\omega(q,p;k)(c) := a(q,p;k)c_qc_p-a(p+k,q-k;k)c_{p+k}c_{q-k} \label{omega}
\end{equation}
where the first sum in the right-hand side of \eqref{DGEDbalanced} arises from $Q_{3,i}(c)+Q_{4,i}(c)$ and the
second term from $Q_{1,i}(c)+Q_{2,i}(c)$.

\medskip

Writing the equations in this form immediately shows that if each of the above reversible chemical reactions
are in equilibrium then all $\omega(q,p;k)(c)=0$ and we are at an equilibrium of the ordinary differential equation system. 
This situation corresponds to the occurrence of the so called microscopic reversibility, and is a physically natural assumption,
translated in mathematical terms in the following \emph{detailed balance}\/ condition: there exists a positive sequence 
$\mathscr{O}=(\mathscr{O}_i),$ with $\mathscr{O}_0=1$, such that $\omega(q,p;k)(\mathscr{O}) = 0, \forall q, p, k.$ Clearly,
under this condition the sequence defined by $\bar{c}_i = \mathscr{O}_jz^j$ will be an equilibrium solution of \eqref{DGED}
if $z>0$ is such that $\bar{c}\in X_{0,1}.$

\medskip

From previous works on cluster dynamics' equations of exchange driven \cite{esenvel} or 
coagulation-fragmentation types \cite{dc15}
with detailed balance assumptions we might expect that, under appropriate conditions, the function
\begin{equation}
V(c) := \sum_{i=0}^\infty c_i\Bigl(\log\frac{c_i}{\mathscr{O}_i} -1\Bigr)\label{V(c)}
\end{equation}
can serve as the foundation from which to construct a Lyapunov function for \eqref{DGED}
in an appropriate topology of the phase space $X_{0,1}$. Clearly, this is the case
when the rate coefficients $a(i, j; k)$ are such that \eqref{DGED} becomes one of the above-mentioned cluster systems.

\medskip

In the case of the generalized exchange-driven system studied in the present paper, computing formally we have 
\[
\frac{d}{dt}V(c(t)) = \sum_{i=0}^\infty \dot{c}_i(t)\log\frac{c_i(t)}{\mathscr{O}_i},
\]
and substituting \eqref{DGEDbalanced}--\eqref{omega} into this expression we formally get
\begin{equation}
\frac{d}{dt}V(c(t)) = \sum_{j=0}^\infty\sum_{i=1}^\infty\sum_{k=1}^i \omega(i,j;k)(c(t))\log\frac{c_{j+k}(t)\mathscr{O}_j}{c_j(t)\mathscr{O}_{j+k}}. \label{dVdt}
\end{equation}

\medskip

Making rigorous these formal computations, studying the sign of \eqref{dVdt}, and hopefully using these tools in the
analysis of the dynamic behaviour of solutions to \eqref{DGED} is clearly the subject of another work.


\begin{thebibliography}{23}




\bibitem{Ball:1990} J. M. Ball and J. Carr, The discrete coagulation-fragmentation 
equations: Existence, uniqueness, and density conservation, J. Stat. Phys., { \bf 61} (1990) 203--234.


\bibitem{Ball:1986} J. M. Ball, J. Carr, and O. Penrose, The Becker-D{\"o}ring cluster 
equations: basic properties and asymptotic behaviour of solutions, Commun. Math. Phys., { \bf 104} (1986) 657--692. 


\bibitem{banakar} Z. Banakar, M.  Tavana, B. Huff, D. Di Caprio, A bank merger 
predictive model using the Smoluchowski stochastic coagulation equation and reverse engineering, 
 Int. J.  Bank Marketing, \textbf{36} (4), (2018) 634--662.

\bibitem{BLLvol1} J. Banasiak, W. Lamb, and Ph. Lauren\c{c}ot, 
Analytic Methods for Coagulation-Fragmentation Models, volume I, 
Monographs and Research Notes in Mathematics, CRC Press, New York, 2019.

\bibitem{BLL} J. Banasiak, W. Lamb, and Ph. Lauren\c{c}ot, 
Analytic Methods for Coagulation-Fragmentation Models, volume II, 
Monographs and Research Notes in Mathematics, CRC Press, New York, 2019.

\bibitem{Ben:2003} E. Ben-Naim and P. L. Krapivsky, Exchange-driven growth,
Phys. Rev. E, {\bf 68} (2003), 031104.

\bibitem{carr92} J. Carr, Asymptotic behaviour of solutions to the coagulation-fragmentation equations. 
I. The strong fragmentation case, Proc. Roy. Soc. Edinburgh Sect. A, \textbf{121} (1992) 231--244.

\bibitem{dc15} F. P. da Costa, Mathematical aspects of
  coagulation-fragmentation equations, in: {\em Mathematics of Energy
and Climate Change}, J.-P. Bourguignon, R. Jeltsch, A. A. Pinto, and
  M. Viana, eds., Springer--Verlag, Cham 2015, pp. 83--162.   
  
\bibitem{degond} P. Degond, J.-G. Liu, R. L. Pego,  Coagulation-fragmentation model for animal group-size statistics, 
 \emph{J. Nonlinear Sci.}, \textbf{27} (2017) 379--424.

 \bibitem{vanDongen} P. G. J. van Dongen and M. H. Ernst,  Kinetics of reversible polymerization, 
 \emph{Journal of Statistical Physics}, \textbf{37} (1984) 301--324.

 \bibitem{Eichenberg:2021}   C. Eichenberg and A. Schlichting,
   Self-similar behavior of the exchange-driven growth model with product kernel,
 Comm. Partial Diff. Eqs., \textbf{46} (2021) 498--546.

\bibitem{esen1} E. Esenturk, 
Mathematical theory of exchange-driven growth, Nonlinearity, {\bf 31} (2018) 3460--3483.

\bibitem{esen2}
E. Esenturk and C. Connaughton, Role of zero clusters in exchange-driven growth with and without input,
Phys. Rev. E, {\bf 101}, (2020) 052134.

\bibitem{esenvel} E. Esenturk and J. Velazquez, Large time behavior of exchange-driven growth,
Discr. Cont. Dyn. Syst., {\bf 41}, 2 (2021) 747--775.

\bibitem{guy} R. D. Guy, A. L. Fogelson, J. P. Keener, Fibrin gel formation in a shear flow,
 Math. Med. Biol., \textbf{24} (1), (2007) 111--130.

\bibitem{HS} P.-F. Hsieh and Y. Sibuya, Basic Theory of Ordinary Differential Equations, 
Universitext, Springer-Verlag, New York, 1999.

\bibitem{KL:2002} J. Ke and Z. Lin, Kinetics of migration-driven aggregation processes with birth and death, Phys.
Rev. E, {\bf 67}, (2002) 031103.

\bibitem{PL:2002} Ph. Lauren\c{c}ot, The discrete coagulation equations with multiple fragmentation,
Proc. Edin. Math. Soc., {\bf 45}, (2002), 67--82.

\bibitem{LR:2002} F. Leyvraz and S. Redner, Scaling theory for migration-driven aggregate growth, Phys. Rev.
Lett.,  {\bf 88}, (2002) 068301.

\bibitem{IKR:1998} S. Ispolatov, P. L. Krapivsky, and S. Redner, Wealth distributions in models of capital exchange,
Eur. J. Phys. B, {\bf 2}, (1998) 267.

\bibitem{mulheran} P .A. Mulheran, Theory of cluster growth on surfaces, in: \emph{Metallic Nanoparticles}, 
J. A. Blackman (ed.), Handbook of Metal Physics, Vol. 5, Elsevier, Amsterdam, 2008, pp. 73--111.

\bibitem{penrose89} O. Penrose, Metastable states for the Becker-D\"oring cluster equations,
Commun. Math. Phys., {\bf 124} (1989) 515--541. 

\bibitem{Pruppacher} H. R Pruppacher and J. D. Klett, Microphysics of Clouds and Precipitation, 2nd Edition,
Atmospheric and Oceanographic Sciences Library, volume 18,
Springer, Dordrecht, 2010.

\bibitem{safronov} V. Safronov, Evolution of the protoplanetary cloud and formation of the earth and the planets,
 Israel Program for Scientific Translations, Jerusalem, 1972.

\bibitem{schl1} A. Schlichting, The exchange-driven growth model: basic properties and longtime behavior,
J. Nonl. Sci., {\bf 30} (2020) 793--830.

\bibitem{Smoluchowski:1917} M. Smoluchowski, Versuch einer mathematischen Theorie der Koagulationskinetik kol-
loider Lösungen, Z. Phys. Chem., {\bf 92} (1917) 129--168.

\end{thebibliography}
\end{document}